\documentclass[10pt]{article}
\usepackage{amsmath,latexsym,amsfonts,amssymb}
\usepackage[latin1]{inputenc}
\newtheorem{theorem}{Theorem}[section]
\newtheorem{lemma}[theorem]{Lemma}
\newtheorem{proposition}[theorem]{Proposition}
\newtheorem{coro}[theorem]{Corollary}
\newtheorem{definition}[theorem]{Definition}
\begin{document}
\newcommand{\pn}{\mathbb P (E)}
\newcommand{\pr}{{\mathbb R\mathbb P}^1}
\newcommand{\cdiff}{C^{1,h}(S^1)\rtimes Di\!f\!f^{h}(S^1)}
\newcommand{\ocdiff}{\Omega^{h}(S^1)\rtimes Di\!f\!f^{h}(S^1)}
\newcommand{\ci}{C^\infty}
\newcommand{\jsr}{J^1(S^1,\mathbb R)}
\newcommand{\pnd}{\mathbb P (E^*)}
\newcommand{\pnu}{\mathbb P (\mathbb R^{n})}
\newcommand{\pnud}{\mathbb P (\mathbb R^{n*})}
\newcommand{\sln}{SL(n,\mathbb R)}
\newcommand{\sld}{SL(2,\mathbb R)}
\newcommand{\inn}{_{n\in\mathbb N}}
\newcommand{\qf}{\mathcal{QF}}
\newcommand{\vm}{\vert\mu\vert}
\newcommand{\Hom}{\rm{Hom}}
\newcommand{\Rep}{\rm{Rep}}
\newcommand{\proof}{{\sc Proof : }}
\newcommand{\rmks}{{\sc Remarks: }}
\newcommand{\rmk}{{\sc Remark: }}
\newcommand{\qed}{{\sc Q.e.d.}}
\newcommand{\grf}{\pi_1 (S)}
\newcommand{\bgrf}{\partial_\infty \pi_1 (S)}
\newcommand{\seq}[1]{
\{#1\}_{m\in\mathbb N}}
\newcommand{\mapping}[4]
{
\left\{
\begin{array}{rcl}
#1 &\rightarrow& #2\\
#3 &\mapsto& #4 
\end{array}
\right.
}
\newcommand{\auteur}{\vskip 2truecm
\centerline{François Labourie} 
\centerline{Topologie et Dynamique} 
\centerline{Université Paris-Sud}
\centerline{F-91405 Orsay (Cedex)}
} 
\title{Cross Ratios,  Surface Groups, $SL(n,\mathbb R)$ and $C^{1,h}(S^1)\rtimes Di\!f\!f^{h}(S^1)$.}
\author{François LABOURIE \thanks{L'auteur remercie l'Institut Universitaire de France.}} 
\maketitle
\begin{abstract}
We  present in this article  relations between cross ratios and representations of surface groups. More specifically, we describe  a connected component of representations into $SL(n,\mathbb R)$ as the set of cross ratios on the boundary at infinity of the group which satisfy some functional relations depending on $n$. We also show that representations into $\cdiff$ can be described as cross ratios. We  exhibit a "character variety" of such representations. We show that this character variety contains all character varieties for $SL(n,\mathbb R)$ as well as  the set of all negatively curved metrics on the surface.
\end{abstract}

\section{Introduction}

In this article, $S$ is a closed surface with genus at least 2, $\grf$ is the fundamental group of $S$. We denote by $\bgrf$ the boundary at infinity of $\grf$. We recall that $\bgrf$ is a one dimensional compact connected Hölder manifold, hence Hölder homeomorphic to the circle $S^1$. The group  $\grf$ acts by Hölder homeomorphisms on $\bgrf$.

We first give a definition, closely related to the beautiful one introduced by Otal in \cite{JPO}.

A {\em (strict) cross ratio} is a $\grf$-invariant  function $b$  defined on  $$
\bgrf^{4*}=\{(x,y,z,t)\in\bgrf^4,\ x\not=t \hbox{ and } y\not= z\}$$
which satisfies some algebraic  rules:
({\em cf.} Equations (\ref{birrules})). 
\begin{eqnarray*}
b(x,y,z,t)&=&b(z,t,x,y) \\
b(x,y,z,t)&=&0\ \ \Leftrightarrow x=y \hbox{ or } z=t\\
b(x,y,z,t)&=&{b(x,y,z,w)}{b(x,w,z,t)}\\
b(x,y,z,t)&=&b(x,y,w,t)b(w,y,z,t)\\
b(x,y,z,t)&=&1 \ \Leftrightarrow x=z \hbox{ or } y=t\
\end{eqnarray*}
Whenever $S$ has a hyperbolic metric, $\bgrf$ is identified with $\pr$. Thus, the usual cross ratio of the projective line gives rise to a cross ratio on $\bgrf$. It follows that every cocompact representation of $\grf$ in $PSL(2,\mathbb R)$ give rises to a cross ratio. Moreover, these  cross ratios  can be characterised as those which satisfies the following extra relation
\begin{eqnarray}
1-b(x,y,z,t)=b(t,y,z,x).\label{crossproj0}
\end{eqnarray}
In this paper, we will in particular generalise this construction when one replaces $SL(2,\mathbb R)$ with $SL(n,\mathbb R)$. We will also comment on  an infinite dimensional version of this construction that appears when one replaces $SL(n,\mathbb R)$ with $\cdiff$. We will also show that the corresponding "character variety" of representations in $\cdiff$ contains a connected component of the character variety of $SL(n,\mathbb R)$ as well the moduli space of negatively curved metrics. We now explain some background and give more precise results in the next two paragraphs.

In \cite{FL3}, we define a {\em $n$-Fuchsian} representation of $\grf$ to be a
representation $\rho$ which may be written as
$\rho=\iota\circ\rho_{0}$, where $\rho_{0}$ is a cocompact
representation of $\grf$ with values in $PSL(2,\mathbb R)$ and $\iota$
is the irreducible representation of $PSL(2,\mathbb R)$ in
$PSL(n,\mathbb R)$.  

In \cite{H}, Hitchin proves the remarkable result that the connected components $\Rep_H(\grf,PSL(n,\mathbb R))$ of the space of reducible representations of $\grf$ in $PSL(n,\mathbb R)$ containing $n$-Fuchsian representations are diffeomorphic to balls. Such a connected component is called a {\em Hitchin} component. It is   denoted by $\Rep_H(\grf,PSL(n,\mathbb R))$. A representation which belongs to a Hitchin component is called a {\em $n$-Hitchin representation}. In other words it is a representation that can be deformed to a $n$-Fuchsian representation.

In \cite{FL3}, we give a geometric description of Hitchin representations, which was later completed by the work of O. Guichard  \cite{Gu} ({\em cf.} Section \ref{hyperconvex}). In particular, we  show that if $\rho$ is a Hitchin representation and $\gamma$ a nontrivial element of $\grf$, then $\rho(\gamma)$ is real split (Theorem 1.5 of \cite{FL3}). 

We finally state another construction about cross ratios: associated to a nontrivial element $\gamma$ of $\grf$ and a cross ratio $b$ is a real number  $l_b(\gamma)$ called the {\em period} of $\gamma$ ({\em cf.} Equation (\ref{period})). For the case of cross ratios associated to hyperbolic metrics, these periods are the length of the associated closed geodesics. 

\subsubsection*{Cross ratios and Hitchin representations: a correspondence}

Our first result describes Hitchin representations in terms of cross ratios, generalising the situation about $PSL(2,\mathbb R)$ which we briefly described in the previous paragraph.

We first introduce more complicated functions built out of cross ratio. For every $p$, let  $\bgrf^p_*$ be the set of pairs of $p+1$-uples
 $(e_0,e_1,\ldots,e_p)$, $(u_0,u_1,\ldots,u_p)$ in $\bgrf$ such that
\begin{eqnarray*}
j>i>0 \implies e_j\not= e_i\not= u_0, u_j\not=u_i\not=e_0.
\end{eqnarray*}
Let $b$ be a cross ratio. Let  $\chi_{b}^p$ be the map  from $\bgrf^p_*$ to $\mathbb R$ defined by
$$
\chi_{b}^p(e,u)=\det_{i,j>0}((b(e_i,u_j,e_0,u_0)).
$$
Our main result is the following. 
\begin{theorem}\label{hcintro}
There exists a bijection between the set of  $n$-Hitchin  representations and the set of cross ratios such that 
\begin{itemize}
\item $\forall e, u$, $\chi_{b}^n(e,u)\not=0$,
\item $\forall e, u $, $\chi_{b}^{n+1}(e,u)= 0$.
\end{itemize}
Furthermore, if $\rho$ is a $n$-Hitchin  representation, $b$ its associated cross ratio, and $\gamma$ a nontrivial element of $\grf$ then the period of $\gamma$ is given by
$$
l_b(\gamma)=\log(\vert\frac{\lambda_{max}(\rho(\gamma))}{\lambda_{min}(\rho(\gamma))}\vert),
$$
where $\lambda_{max}(\rho(\gamma))$ and $\lambda_{min}(\rho(\gamma))$ are respectively the  eigenvalues  of respectively maximum and minimum absolute values of the element $\rho(\gamma)$.
\end{theorem}

For $n=2$, it turns out that the functional relations described in this Theorem amounts to Relation \ref{crossproj0}

The {\em  limit curve}, drawn in $\pnu$ and described in Paragraph \ref{limitcurve}, is the link between cross ratio and representations.

\subsubsection*{A "Character variety" containing all Hitchin representations}
 Since Hitchin representations are irreducible, the natural  embedding of $PSL(n,\mathbb R)$  in $PSL(n+1,\mathbb R)$ does not give rise to an embedding of the corresponding Hitchin components. Therefore there is no natural algebraic way -  by an injective limit procedure say - to build a limit when $n$ goes to infinity of Hitchin components. However,  it follows  from the previous paragraph that all Hitchin components sit in the same "moduli space": the space of all cross ratios.

We now explain the second construction of this article:  {\em all Hitchin components lie in a "character variety" of $\grf$ into an infinite dimensional group $G$}. More precisely, let $C^{1,h}(S^1)$ be the vector space of $C^1$-functions with Hölder derivatives on the circle, and let $Di\!f\!f^h(S^1)$ be the group of  $C^1$-diffeomorphisms with Hölder derivatives of the circle. We observe that  $Di\!f\!f^h(S^1)$ acts naturally on $C^{1,h}(S^1)$. Let 
$$G=\cdiff.$$   
In Paragraph \ref{actionj1}, we explain that this group $G$ has a natural action by Hölder homeomorphisms on $J^1(S^1)$,  the space of $1$-jets of functions  on the circle.

Our first result is to build a "character variety" for homomorphisms of $\grf$ in $G$. Namely, we define in Paragraph \ref{goodhomo} {\em $\infty$-Hitchin homomorphisms }of $\grf$ in $G$. For such a $\infty$-Hitchin homomorphism $\rho$, the quotient  $J^1(S^1)/\rho(\grf)$ is compact. Moreover, in Paragraph \ref{specthomh}  we associate a real number $l_\rho(\gamma)$ called the {\em  $\rho$-length} of $\gamma$ to every nontrivial element $\gamma$ of $\grf$, and every $\infty$-Hitchin representation $\rho$. The marked collection of $\rho$-lengths is the {\em spectrum} of $\rho$. In Paragraph \ref{crosshomh}, we also associate to every $\infty$-Hitchin homomorphism a cross ratio. 

We denote by $\Hom_H$ the set of all $\infty$-Hitchin representations. Let $Z(G)$ be the center of $G$. Our first result describe the action on $G$ on $\Hom_H$.

\begin{theorem}\label{homhintro}
The set $\Hom_H$ is an open set in the set $\Hom(\grf,G)$. Moreover, $G/Z(G)$ acts properly on $\Hom_H$ and the quotient $\Hom(\grf,G)/G$ is Haussdorf. Two representations with the same spectrum and the same cross ratio are conjugated.
\end{theorem}

The two properties of $\infty$-Hitchin representations stated in the previous theorem show they are good candidates to be avatars of reducible representations. We denote by $\Rep_H$ the {\em character variety} of $\infty$-Hitchin representations $\Hom_H/G$. The following result relates Hitchin components to this  character variety.

\begin{theorem}\label{hitchcrossintro}
There exists a continuous injective  map 
$$
\psi : \Rep_{H}(\grf,SL(n,\mathbb R))\rightarrow \Rep_{H},
$$
such that, if $\rho\in \Rep_{H}(\grf,SL(n,\mathbb R))$ then
\begin{itemize}
\item for any $\gamma$ in $\grf$, we have 
\begin{equation}
l_{\psi(\rho)}(\gamma)=\log(\vert\frac{\lambda_{max}(\rho(\gamma)}{\lambda_{min}(\rho(\gamma)}\vert).
\end{equation}
Here, $l_{\psi(\rho)}(\gamma)$ is the $\psi(\rho)$-length of $\gamma$, and $\lambda_{max} (a)$ (resp. $\lambda_{min}(a)$) denote the maximum (resp. minimum) real eigenvalue of the endomorphism $a$ in absolute value.
\item The cross ratio associated to $\rho$ and $\psi(\rho)$ coincide.
\end{itemize}
\end{theorem}
In some sense, this result  says that $\cdiff$ is a version of $SL(\infty,\mathbb R)$.

We finally state in this introduction another result that explains that our character variety contains yet another interesting space.

\begin{theorem}\label{negcurvintro}
Let $\mathcal M$ be the space of negatively curved metrics on the surface $S$. There exists a continuous injective map $\psi$ from $\mathcal M$ to $\Rep_{H}$. Furthermore, $\psi$ preserves the length: for any $\gamma$ in $\grf$
$$
l_{g}(\gamma)=l_{\psi(g)}(\gamma).
$$
Here $l_{g}(\gamma)$ is the length of the closed geodesic for $g$ freely homotopic to $\gamma$, and $l_{\psi(g)}$ is the $\psi(g)$
-length of $\gamma$.
Finally, $\psi(g_{0})=\psi(g_{1})$, if and only if there exists a diffeomorphism $F$of $S$,  homotopic to the identity, such that $F^*(g_{0})=g_{1}$. 
\end{theorem}

Both results are consequences of a general conjugation result: Theorem \ref{conjugtheo}.

We finish this general introduction by stating a question about our construction: can one characterise 
$$
F={\overline{ \bigcup_n \Rep_H(\grf,PSL(n,\mathbb R))}}
$$
in $\Rep_H$ ? For instance does $F$ contains $\mathcal M$ ?

\subsubsection*{Structure of the article}
We now describe shortly the content of this article.
\begin{itemize}
\item[\bf{\ref{hyperconvex}}.]{\em Curves and hyperconvex representations.} We recall results of \cite{FL3} and explain how Hitchin representations are related to special curves in projectives spaces.
\item[\bf{\ref{cross ratiodef}}.]{\em Cross ratio, definitions and first properties.} We give the precise definition of a cross ratio and of the related quantities (periods and triple ratio).
\item[\bf{\ref{cross ratio}}.]{\em Examples of cross ratio.} We explain various constructions of cross ratio: the classical cross ratio on the projective line, cross ratio associated to curves in projective spaces, dynamical cross ratio and the original construction of J.-P. Otal for negatively curved metrics.
\item[\bf{\ref{hitchincross}}.]{\em Hitchin representations and cross ratios.} We prove Theorem \ref{hcintro}.
\item[\bf{\ref{jetspacesect}}.]{\em The jet space $J^1(S^1,\mathbb R)$.} We describe the geometry of this jet space and the action of $\cdiff$ on it.
\item[\bf{\ref{homosect}}.]{\em Homomorphisms of $\grf$ in $\cdiff$.} We give the precise definitions of Anosov and $\infty$-Hitchin representations of $\grf$ and $\cdiff$, and of the {\em spectrum}  and cross ratio associated to these representations. We prove a refinement of Theorem \ref{homhintro}.
\item[\bf{ \ref{conjugtheosect}}.]{\em A Conjugation Theorem.} We explain Theorem \ref{conjugtheo} which shows how to build  $\infty$-Hitchin representations of $\grf$ in $\cdiff$
\item[\bf{ \ref{negcurvsect}}.]{\em Negatively curved metrics.} We prove Theorem \ref{negcurvintro}.
\item[\bf{\ref{hichsect}}.]{\em Hitchin component.} We prove Theorem \ref{hitchcrossintro}.
\item[\bf{\ref{applam}}.]{\em Appendix A: Filtrated Spaces, Holonomy.} We prove results about  foliations by affine spaces.
\item[\bf{\ref{sympnat}}.]{\em Appendix B: the symplectic nature of cross ratio.} We explain the relation of cross ratios are related and symplectic geometry, and give constructions of a whole family of cross ratio related to Hitchin representations.
\end{itemize}
\section{Curves and hyperconvex representations}\label{hyperconvex}

We recall results and definitions from \cite{FL3}.

\subsection{Hyperconvex representations}
\subsubsection{Fuchsian representations}
A {\em $n$-Fuchsian} representation of $\grf$ is  a
representation $\rho$ which may be written as
$\rho=\iota\circ\rho_{0}$, where $\rho_{0}$ is a cocompact
representation with values in $PSL(2,\mathbb R)$ and $\iota$
is the irreducible representation of $PSL(2,\mathbb R)$ in
$PSL(n,\mathbb R)$. 

\subsubsection{Hyperconvex curves} \label{limitcurve}
A continuous curve $\xi$  with values in
$\mathbb P(\mathbb R^{n})$ is  {\em hyperconvex}
if, for any distinct points $(x_{1},\ldots,x_{n})$, the
following sum is direct
$$
\xi(x_{1})+\ldots+\xi(x_{n}).
$$
We say a representation $\rho$ of $\grf$ is {\em$n$-hyperconvex}, if there exists a $\rho$-equivariant hyperconvex curve from $\partial_\infty\grf$ in $\mathbb  P(\mathbb R^n)$. Actually, such a curve is unique and is called  the {\em limit curve} of the representation.
We say a representation is Hitchin if  may be deformed in a $n$-Fuchsian representation.

In \cite{FL3} we prove the following result.
\vskip 0.5truecm
\noindent{\bf Theorem. }{\em Let $\rho$ be a Hitchin representation.  Then $\rho$ is hyperconvex.  Each such representation is discrete, faithful. Finally, for each
  $\gamma$ in $\grf$ different from the identity,
  $\rho(\gamma)$ is real split with distinct eigenvalues}
\vskip 0.5truecm
We explain later a refinement of this result (Theorem  \ref{mainanosov}). We observe that the Veronese embedding is a hyperconvex curve equivariant under all Fuchsian representations. Therefore, a Fuchsian representation is indeed hyperconvex.

According to the previous result, many representations, at least those in Hitchin components, are hyperconvex.

Conversely, completing our work, O. Guichard \cite{Gu} has shown the following result

\begin{theorem}
{\sc [Guichard]} Every hyperconvex representation is Hitchin.
\end{theorem}

\subsection{Frenet curves}
We say a hyperconvex curve $\xi$ is a {\em Frenet
  curve}, if there exists a family of maps
$(\xi^{1},\xi^{2},\ldots,\xi^{n-1})$ , called the {\em osculating flag},
such that 
\begin{itemize}
\item $\xi^{p}$
 takes values in the Grassmannian of $p$-planes,
 \item 
$\forall x,\ \ \xi^{p}(x)\subset\xi^{p+1}(x)$
\item $\xi=\xi^{1}$,
\item if $(n_1,\ldots,n_l)$ are positive integers such that $\sum_{i=1}^{i=l}n_i\leq n, $ if $(x_{1},\ldots,x_{l})$ are
  distinct points, then the following sum is direct
\begin{eqnarray}
\xi^{n_i}(x_i)+\ldots+\xi^{n_{l}}(x_{l});
\end{eqnarray}
\item finally, for every $x$, let $p=n_{1}+\ldots+n_{l}$,
  then
\begin{eqnarray}
\lim_{(y_1,\ldots,y_l)\rightarrow x, y_i
{\hbox{\tiny all distinct}}}
(\bigoplus_{i=1}^{i=l}\xi^{n_i}(y_i))=\xi^{p}(x).
\end{eqnarray}
\end{itemize}
We call $\xi_{n-1}$ the {\em osculating hyperplane}. 
We observe that for a Frenet hyperconvex curve, $\xi^{1}$
completely determines $\xi^{p}$. Moreover, if $\xi^{1}$ is
$C^{\infty}$, then $\xi^{p}(x)$ is completely generated by
the derivatives at $x$ of $\xi^{1}$ up to order $p-1$.
However, in general, a Frenet hyperconvex curve has no
reason to be $C^{\infty}$ although its image is obviously a
$C^{1}$-submanifold.

\subsection{Hyperconvex representations and Frenet curves}\label{crosscurve}

We list here several properties of hyperconvex representations proved in \cite{FL3}.

\begin{theorem}\label{mainanosov}
 Let $\rho$ be an hyperconvex representation of $\grf$ in $SL(E)$, with limit curve $\xi$.    \begin{enumerate}
  \item
  Then for each
  $\gamma$ in $\grf$ different from the identity,
  $\rho(\gamma)$ is purely real split with distinct eigenvalues. 
  \item Furthermore, $\xi$ is a hyperconvex Frenet curve. 
  \item Let $\xi^*$ be its osculating hyperplane, then $\xi^*$ is hyperconvex. 
  \item The osculating flag is Hölder.
  \item  Finally, if $\gamma^+$ is the attracting fixed point of $\gamma$ in $\partial_\infty\grf$, then $\xi(\gamma^+)$, (resp.
 $\xi^*(\gamma^+)$) is the unique attracting fixed point of $\rho(\gamma)$ in $\pn$ (resp. $\pnd$)
 \end{enumerate}
\end{theorem}

\section{Cross ratio, definitions and first properties}\label{cross ratiodef}

\subsection{Cross ratio}
Let $S$ be the circle.  Let  
$$S^{4*}=\{(x,y,z,t)\in S^4\ x\not=t , \hbox{ and } y\not= z\}.$$
A {\em cross ratio} on $S$ is a  Hölder function $b$ on $S^{4*}$ with values in $\mathbb R$ which satisfies the following rules
\label{birrules}
\begin{eqnarray}
b(x,y,z,t)&=&b(z,t,x,y) \label{bir100}\\
b(x,y,z,t)&=&0\ \ \Leftrightarrow x=y \hbox{ or } z=t\\
b(x,y,z,t)&=&{b(x,y,z,w)}{b(x,w,z,t)}\label{bir11}\\
b(x,y,z,t)&=&b(x,y,w,t)b(w,y,z,t)\label{bir11bis}
\end{eqnarray}
If $S=\bgrf$, we assume furthermore that $b$ is invariant under the diagonal action of $\grf$:
\begin{eqnarray}
\forall \gamma\in\grf,\ \  b(\gamma x,\gamma y,\gamma z,\gamma t)=b(x,y,z,t)
\end{eqnarray}
Furthermore, we say a cross ratio on $S$ is {\em strict} if
\begin{eqnarray}
b(x,y,z,t)&=&1 \ \Leftrightarrow x=z \hbox{ or } y=t\label{bir12}.
\end{eqnarray}
The classical cross ratio on $\pr$ is an example of a strict cross ratio. It is a well known fact that the classical cross ratio can be characterised as the unique cross ratio satisfying an extra functional rule ({\em cf} Proposition \ref{projline}). We give more examples and constructions in Section \ref{cross ratio}.
\vskip 1 truecm
\noindent\rmk

 {\em The definition given above does not coincide with the usual definition given for instance in \cite{Ham}, \cite{led} (even after taking an exponential): indeed, first we require $b(z,t,x,y)=b(x,y,z,t)$, (which amounts to time reversibility), more importantly we do not require $b(x,y,z,t)=b(y,x,t,z)$.
However we may observe that if $b(x,y,z,t)$ is a cross ratio with our definition, so is $b^*(x,y,z,t)=b(y,x,t,z)$, and finally $bb^*$ is a cross ratio according to the classical definitions quoted above. We explain Otal's construction and the relation with negatively curved metrics in Section \ref{dynot}}

\subsection{Periods}
Let $b$ be a cross ratio $b$ and $\gamma$ be a nontrivial element in $\grf$. The {\em period} $l_b(\gamma)$ is defined as follows.  Let  $\gamma^{+}$ (resp. $\gamma^{-}$) the attracting (resp.  repelling) fixed point of $\gamma$ on $\partial_\infty\pi_1 (S)$.  Let $y$ be an element of $\partial_\infty \pi_1 (S)$. Let's define
\begin{eqnarray}
l_b(\gamma,y)=\log\vert b(\gamma^-,\gamma y,\gamma^+,y)\vert.\label{period}
\end{eqnarray}
It is immediate to check that $l_b(\gamma)=l_b(\gamma,y)$ does not depend on $y$. Moreover, by Equation( \ref{bir100}), $l_b(\gamma)=l_b(\gamma^{-1})$.

\subsection{Triple ratios}\label{triratio}
We introduce a new feature associated to a cross ratio: for every quadruple of pair-wise distinct points $(x,y,z,t)$, one easily checks that the expression
$$
b(x,y,z,t)b(z,x,y,t)b(y,z,x,t),
$$
is independent of the choice of $t$. We call such a function a {\em triple ratio}. Indeed, in some cases, it is related to the triple ratios introduced by A. Goncharov in \cite{G}. It turns out, although we do not use this remark, that a triple ratio satisfies the (multiplicative) cocycle identity and hence defines a bounded cohomology class in $H^2_b(\Gamma)$.

\section{Examples of  cross ratio}\label{cross ratio}
\subsection{Cross ratio on the projective line}
Let $E$ be a vector space with $\dim(E)=2$. We recall that the "classical" cross ratio is defined on $\mathbb P (E)$, identified with $\mathbb R\cup \{\infty\}$ using  projective coordinates, by:
$$
b(x,y,z,t)=\frac{(x-y)(z-t)}{(x-t)(z-y)}.
$$
It is easy to check that this "classical" cross ratio is a strict cross ratio.
The "classical" cross ratio on the projective line satisfies  the rules \ref{birrules} as well as the following  extra rule :
\begin{eqnarray}
1-b(f,v,e,u)=b(u,v,e,f)\label{cr0}.
\end{eqnarray}
We will later on  explain the well known fact that this extra relation completely characterises the classical cross ratio. Furthermore, it turns out that this (simple) relation is equivalent to the following more sophisticated one, which we may generalise to higher dimensions
\begin{proposition}\label{projline}
For a cross ratio $b$, Relation (\ref{cr0}) is  equivalent to,
 \begin{eqnarray}
(b(f,v,e,u)-1)(b(g,w,e,u)-1)=(b(f,w,e,u)-1)(b(g,v,e,u)-1)\label{cr1}.
\end{eqnarray}
Furthermore, if  a cross ratio $b$ satisfies Relation (\ref{cr1}),  there exists an embedding $f$, unique up to left composition to projective transformations, of $S$ in $\mathbb{RP}^1$, such that in projective coordinates,
$$
b(x,y,z,t)=\frac{(f(x)-f(y))(f(z)-f(t))}{(f(x)-f(t))(f(z)-f(y)}.
$$
\end{proposition}
\proof First, if $b$ satisfies Relation (\ref{cr0}), then thanks to Relation (\ref{bir11}), it satisfies (\ref{cr1}).
Conversely, assume it satisfies (\ref{cr1}), then
\begin{eqnarray}
b(f,v,e,u)&=&\frac{(b(f,w,e,u)-1)(b(g,v,e,u)-1)}{b(g,w,e,u)-1}+1.
\end{eqnarray}
Setting $g=v$, we obtain
\begin{eqnarray}
b(f,v,e,u)
&=&\frac{1-b(f,w,e,u)}{b(v,w,e,u)-1}+1\\
&=&\frac{b(v,w,e,u)-b(f,w,e,u)}{b(v,w,e,u)-1}.\label{bir122}
\end{eqnarray}
We have
\begin{eqnarray}
b(f,v,k,z)
&=&\frac{b(f,v,e,z)}{b(k,v,e,z)} {\hbox{ by (\ref{bir11bis})}}
\\
&=&\bigg(\frac{b(f,v,e,u)}{b(f,z,e,u)}\bigg)\bigg(\frac{b(k,z,e,u)}{b(k,v,e,u)}\bigg){\hbox{ by  (\ref{bir11})}}\\
&=&\frac{b(f,v,e,u)b(k,z,e,u)}{b(k,v,e,u)b(f,z,e,u)}\label{bir222}.
\end{eqnarray}
Finally, applying Relation (\ref{bir122}) to the four left terms of Equation (\ref{bir222}) we obtain
$$
b(f,v,k,z)=\frac{(b(v,w,e,u)-b(f,w,e,u))(b(z,w,e,u)-b(k,w,e,u))}{(b(v,w,e,u)-b(k,w,e,u))(b(z,w,e,u)-b(f,w,e,u))}.
$$
The final statement follows when we set
$$
f=f_{(w,e,u)}\mapping{S}{\mathbb R\cup\{\infty\}}{x}{b(x,w,e,u)}.
$$
\qed

\subsection{Cross ratios and curves in projective spaces}
We now extend the previous discussion. 
Let $E$ be an $n$-dimensional vector space.
Let $\xi$ and $\xi^*$ be two curves from  $S$ to $\mathbb P(E)$ and $\mathbb P(E^*)$ respectively. We assume furthermore 
\begin{eqnarray}
\langle\xi^*(z),\xi(y)\rangle=0 \Leftrightarrow z=y \label{condicurve}.
\end{eqnarray}
This in particular true for an equivariant hyperconvex curve $\xi$ and its osculating hyperplane $\xi^*$. 

For every $x$, we choose an arbitrary nonzero vector $\hat\xi(x)$ (resp. $\hat\xi^*(x)$) in the line $\xi(x)$ (resp. $\xi^*(x)$).

We define the cross ratio associated to this pair of curves by
$$
b_{\xi,\xi^*}(x,y,z,t)=\frac{\langle\hat\xi(x),\hat\xi^*(y)\rangle\langle\langle\hat\xi(z),\hat\xi^*(t)\rangle}{\langle\hat\xi(z),\hat\xi^*(y)\rangle\langle\hat\xi(x),\hat\xi^*(t)\rangle}.
$$
It is easy to check that
\begin{itemize}
\item this definition does not depend on the choice of $\hat\xi$ and $\hat\xi^*$, 

\item $b_{\xi,\xi^*}$ satisfies the axioms of a formal cross ratio.
\item Let $V=\xi(x)\oplus\xi(z)$. Let $\eta(m)=\xi^*(m)\cap V$. Let $b_V$ be the classical cross ratio on $\mathbb P(V)$, then
$$
b_{\xi,\xi^*}(x,y,z,t)=b_V(\xi(x),\eta(y),\xi(z),\eta(t)).
$$  
\item It follows that $b_{\xi,\xi^*}$ is strict if furthermore, for all quadruple of pairwise distinct points $(x,y,z,t)$,
$$
Ker(\xi^*(z)) \cap \big(\xi(x)\oplus \xi(y)\big)\not=Ker(\xi^*(t)) \cap 
\big(\xi(x)\oplus \xi(y)\big) \label{injcross}.
$$
\end{itemize}

Finally, we have
\begin{lemma} Let $(\xi,\xi^*)$ and $(\eta,\eta^*)$ be two pairs of curves satisfying Condition (\ref{condicurve}). Assume that $\xi^*$ and $\eta*$ are hyperconvex. Assume furthermore that
$$
b_{\eta,\eta^*}=b_{\xi,\xi^*}.$$ 
Then there exists a linear map $A$ such that
$\xi=A\eta$. \label{unixi}
\end{lemma}
\proof We assume the hypothesis of the theorem. Let $(x_0,x_1,\ldots,x_n)$ be a tuple of $n+1$ pair-wise distinct points of $S$. We choose a vector $z_0$ in $\xi^*(x_0)$. Let $U=(u_1,\ldots,u_n)$ be the basis of $E^*$ such that
$$
u_i\in \xi^*(x_i) \hbox{ and } \langle  z_0,u_i \rangle=1.
$$
The projective coordinates of $\xi (y)$ in the dual basis are
\begin{eqnarray*}
[\ldots:\langle \xi(y),u_i\rangle:\ldots]&=&[ \ldots: \frac{\langle \xi(y),u_i \rangle}{\langle \xi(y),u_1 \rangle}:\ldots ]\cr
&=&[ \ldots:\frac{\langle \xi(y),u_i \rangle}{\langle \xi(y),u_1\rangle} \frac{\langle z_0,u_1 \rangle}{\langle z_0,u_i\rangle}:\ldots ]\cr
&=&[\ldots: b_{\xi,\xi^*}(y,x_i,x_0,x_1):\ldots].
\end{eqnarray*}
Symmetrically, we choose  a vector $y_0$ in $\eta^*(x_0)$, and let $V=(v_1,\ldots,v_n)$ be the basis of $E^*$ such that
$$
v_i\in \eta^*(x_i) \hbox{ and } \langle  y_0,v_i \rangle=1.
$$
It follows that if $A$ is the linear map that send the dual basis of $V$ to the dual basis of $U$, then $\xi=A.\eta$.\qed

We state the following Proposition.

\begin{proposition}
Let $\rho$ be a  a hyperconvex representation of $\grf$ in $SL(n,\mathbb R)$, with limit curve $\xi$. Let  $\xi^*$ be its osculating hyperplane (cf. Theorem \ref{mainanosov}). Then $b_\rho=b_{\xi,\xi^*}$ is a  strict cross ratio defined on $\bgrf$ and its periods are
$$
l_b(\gamma)=\log(\vert\frac{\lambda_{max}(\rho(\gamma)}{\lambda_{min}(\rho(\gamma)}\vert),
$$
where $\lambda_{max}(\rho(\gamma)$ and $\lambda_{min}(\rho(\gamma)$ are respectively the  eigenvalues  of respectively maximum and minimum absolute value  of the purely loxodromic element $\rho(\gamma)$.
\end{proposition}
\proof We already now that $b_\rho$ is a formal cross ratio.  We will prove  in Proposition \ref{b=1} that $b_\rho$ is strict. We now compute the periods. By Theorem \ref{mainanosov}, if $\gamma^+$ is the attracting fixed point of $\gamma$ in $\partial_\infty\grf$, then $\xi(\gamma^+)$, (resp.
 $\xi^*(\gamma^-)$) is the unique attracting (resp. repelling) fixed point of $\rho(\gamma)$ in $\pn$. In particular 
\begin{eqnarray*}
\rho(\gamma)\hat\xi(\gamma^+)&=&\lambda_{max}\hat\xi(\gamma^+)\\
\rho(\gamma)\hat\xi(\gamma^-)&=&\lambda_{min}\hat\xi(\gamma^-)\\
\end{eqnarray*}
Therefore,
\begin{eqnarray*}
l_b(\gamma)&=&\log\vert b(\gamma^{-},y,\gamma^{+},\gamma^{-1} y)\vert\\
&=&\log\big\vert\frac{\langle\hat\xi(\gamma^{-}),\hat\xi^*(y)\rangle \langle\hat\xi(\gamma^{+}),\hat\xi^*(\gamma^{-1}y)\rangle}{\langle\hat\xi(\gamma^{-}),\hat\xi^*(\gamma^{-1} y)\rangle \langle\hat\xi(\gamma^{+}),\hat\xi^*(y)\rangle}\big\vert\\
&=&\log\big\vert\frac{\langle\hat\xi(\gamma^{-}),\hat\xi^*(y)\rangle \langle\hat\xi(\gamma^{+}),\rho(\gamma)^*\hat\xi^*(y)\rangle}{\langle\hat\xi(\gamma^{-}),\rho(\gamma)^*\hat\xi^*( y)\rangle \langle\hat\xi(\gamma^{+}),\hat\xi^*(y)\rangle}\big\vert\\
&=&\log\big\vert\frac{\langle\hat\xi(\gamma^{-}),\hat\xi^*(y)\rangle \langle\rho(\gamma)\hat\xi(\gamma^{+}),\hat\xi^*(y)\rangle}{\langle\rho(\gamma)\hat\xi(\gamma^{-}),\hat\xi^*( y)\rangle \langle\hat\xi(\gamma^{+}),\hat\xi^*(y)\rangle}\big\vert\\
&=&\log\big\vert\frac{\lambda_{max}}{\lambda_{min}}\big\vert
\end{eqnarray*}
\qed
\vskip 1truecm
\rmk
Let $\rho^*$ be the contragredient representation of $\rho$ defined by
$$
\rho^*(\gamma)=(\rho(\gamma^{-1}))^*.
$$
Let $b$ be a cross ratio. We define $b^*$ by  $b^*(x,y,z,t)=b(y,x,t,z)$.  Then $(b_\rho)^*=b_\rho^*$ and $b_\rho$ and $b^*_\rho$ have the same periods.
\subsection{Dynamical cross ratio and other examples}\label{dynot}

We recall Otal's construction of cross ratio in the case of negatively curved metrics on surfaces \cite{JPO}. Let $S$ be equipped with a negatively curved metric. We lift this metric to the universal cover $\tilde S$ of $S$. Let $(a_1,a_2,a_3,a_4)$ be a quadruple of points on $\bgrf=\partial_\infty\tilde S$. Let $c_{ij}$ be the geodesic from $a_i$ to $a_j$.  We choose nonintersecting horoballs
$H_i$ "centred" at each point $a_i$. Let $l_{ij}$ be the length of the following geodesic arc 
$$
c_{ij}\cap(\tilde S\setminus (H_i\cup H_j)).
$$
{\em Otal's cross ratio} is
$$
O(a_1,a_2,a_3,a_4)=l_{12}-l_{23}+l_{34}-l_{41}.
$$
Let's make the link with our definition: the cross ratio  of a {\it cyclically oriented 4-uple} of distinct points is just the exponential of Otal's cross ratio
$$
b(a_1, a_2,a_3,a_4)=e^{l_{12}-l_{23}+l_{34}-l_{41}}.
$$
For noncyclically oriented 4-uple, we introduce a sign compatible with  the sign of the usual cross ratio.

With this definition the cross ratio agrees with the cross ratio defined on $\pr$ in the case of  hyperbolic surfaces.
\vskip 1 truecm
\rmks
\begin{enumerate}
\item Actually, this construction can be extended to Anosov flows on unit tangent bundle of surfaces. Conversely, it is very easy to show in the case of surfaces, that a cross ratio comes from a flow. Indeed, given a cross ratio $b$ on $\partial_\infty\grf^{4*}$, and any real number $t$, we define a flow map $\phi_t$ from $\partial_\infty\grf^{3*}$ to itself by
$$
\phi_t(x_-,x_0,x_+)=(x_-,x_t,x_+),
$$
where
$$
b(x_+,x_0,x_-,x_t)=e^t.
$$
Then, Equation (\ref{bir11}) exactly says that $t\mapsto \phi_t$ is a one parameter group.
\item The formal rules \ref{birrules} of a cross ratio are easily satisfied 
\item The associated period of an element $\gamma$ is the length of the closed orbit in the free homotopy class.
\item the cross ratio of a negatively curved metric satisfies an extra symmetry: $$
b(x,y,z,t)=b(y,x,t,z).
$$
This symmetry is required in Otal's definition but not in ours. We pay this generalisation by the fact a cross ratio is not uniquely determined by its periods.
\item Of course it is well known that  this construction is yet another instance of the "symplectic construction" explained in Section \ref{sympnat}. Indeed the space of geodesics of the universal cover of a negatively curved manifold is equipped with a symplectic structure which comes from the symplectic reduction of the tangent space. Furthermore, this space is identified with
the space pair of distinct points of the boundary at infinity. Hence, it inherits a real polarisation and one checks that the cross ratio defined as in Section \ref{sympnat} coincides with the dynamical one.
\end{enumerate} 

For a complete description of various aspects of dynamical cross ratios, one is advised to read François Ledrappier's presentation \cite{led}.

\section{Hitchin representations and cross ratios} \label{hitchincross}
We  aim to generalise Proposition
\ref{projline}  in higher dimensions.

Let $b$ be a cross ratio on a set $S$. We now introduce more  sophisticated functions. 
Let $e=\{e_0,e_1 \ldots e_n\}$ be $n+1$ points on $S$, and $u=\{u_0,u_1,\ldots , u_n\}$ $n+1$ other points. Assume that,
$$
\forall i, e_i\not=u_0,\ \ \forall i, u_i\not=e_0.
$$ 
Let $B(e;u)$ be the $n\times n$-matrix whose coefficients are 
$$
b_{ij}=b(e_i,u_j,e_0,u_0) \hbox{ with } 1\leq i\leq n, \ \ 1\leq j\leq n.
$$
We set
\begin{displaymath}
\chi_b^n(e;u)= \det(B(e;u)).
\end{displaymath}
When the context makes it obvious, we omit the subscript $b$.
We now prove.
\begin{theorem}\label{maincross ratio}
Let  $\xi$ and $\xi_*$ be two hyperconvex curves from $S^1$ to $\pn$ and $\pnd$ such that 
$$
x=y\Leftrightarrow\langle\xi_*(x),\xi(y)\rangle=0.
$$ 
Let $b=b_{\xi,\xi^*}$ be the associated cross ratio.
Then, for all pairs of $p$-uples $(e_0,e_1,\ldots,e_p)$, $(u_0,u_1,\ldots,u_p)$ with
\begin{eqnarray*}
j>i>0 \implies e_j\not= e_i\not= u_0, u_j\not=u_i\not=e_0,
\end{eqnarray*}
we have , 
\begin{enumerate}
\item if $p\leq n$, $\chi_b^p(e,u)\not=0$,
\item if $p >n$, $\chi_b^p(e,u)= 0$. 
\end{enumerate}
Conversely, let $b$ be  a continuous strict cross ratio on $S^1$. Assume that the cross ratio satisfies the above two conditions, then there exist two hyperconvex curves $\xi$ and $\xi^*$ with values in $\pn$ and $\pnd$ respectively, unique up to projective transformations, such that
$b=b_{\xi,\xi^*}$
and
$
x=y\Leftrightarrow\langle\xi_*(x),\xi(y)\rangle=0.
$ 
\end{theorem}
The proof of the Theorem is given in Paragraph \ref{maincross ratio}.
\vskip 1 truecm
\rmks
\begin{itemize}
\item We shall very soon prove that for all cross ratios, the nullity of  $\chi^p(e,u)$ for a given $e$ and $u$  does not depend on $e_0$ and $u_0$,
\item
It is interesting to see what happens when $n=2$. In this case, using the previous remark, let's take $u_0=u_1=u$, $e_0=e_1=e$.
Then 
\begin{eqnarray*}
\chi^2(e,e,f;u,u,v)&=&
\left|
\begin{array}{cc}
1 & 1\\
1 &b(f,v,e,u)
\end{array}
\right|
\\
&=& b(f,v,e,u)-1\not=0.
\end{eqnarray*}
Also,
\begin{eqnarray*}
& &\chi^3(e,e,f,g;u,u,v,w)\\
&=&
\left|
\begin{array}{ccc}
1 & 1& 1\\
1 &b(f,v,e,u)&b(g,v,e,u)\\
1&b(f,w,e,u)&b(g,w,e,u)
\end{array}
\right| \\
&=&
\left|
\begin{array}{ccc}
1 & 0& 0\\
1 &b(f,v,e,u)-1&b(g,v,e,u)-1\\
1&b(f,w,e,u)-1&b(g,w,e,u)-1
\end{array}
\right| \\
&=&(b(f,v,e,u)-1)(b(g,w,e,u)-1)-(b(f,w,e,u)-1)(b(g,v,e,u)-1).
\end{eqnarray*}
Therefore, the condition that $\xi^p=0$ extends Equation (\ref{cr1}). 
\end{itemize}

\subsection{The expression $\chi$}
We first prove the following result of independent interest.
\begin{proposition}
The nullity of  $\chi^n(e,u)$ does not depend on $e_0$ and $u_0$.
\end{proposition}
\proof
From the formal rules of computation of cross ratios, we have
\begin{eqnarray*}
b(e_i,u_j,e_0,u_0)b(e_i,u_0,e_0,v_0)=b(e_i,u_j,e_0,v_0)=b(e_i,u_j,f_0,v_0)b(f_0,u_j,e_0,v_0)
\end{eqnarray*}
Therefore
\begin{eqnarray*}
& &\chi(e_0,e_1,\ldots,e_n;u_0,u_1,\ldots,u_n)\\
&=&\Big(\prod_{i,j}\frac{b(f_0,u_j,e_0,v_0)}{b(e_i,v_0,e_0,u_0)}\Big)\chi(f_0,e_1,\ldots,e_n;v_0,u_1,\ldots,u_n).
\end{eqnarray*}
Since we have assumed
$$
f_0\not=v_0\not=e_i, u_j\not=e_0\not=f_0,
$$
The proposition immediately follows. 
\qed

\subsection{Proof of Theorem \ref{maincross ratio}}
\subsubsection{Cross ratio associated to curves}
Let $\xi$ and $\xi^*$ be two  hyperconvex curves with values in $\pn$ and $\pnd$ such that
$$
x=y\Leftrightarrow\langle\xi_*(x),\xi(y)\rangle=0.
$$ 
Let  $e_0$ and $u_0$ be two  distinct points of $S^1$. Let $E_0$ be a nonzero vector  in $\xi(e_0)$ . Let  $U_0$ be a covector in $\xi^*(u_0)$ such that $\langle U_0,E_0 \rangle=1$.
We now lift the curves $\xi$ and $\xi^*$ with values in $\pn$ and $\pnd$ to continuous curves $\hat \xi$ and $\hat\xi^*$ from $S^1\setminus\{e_0,u_0\}$ to $E$ and $E^*$, such that
\begin{eqnarray*}
\langle\hat\xi^*(v),E_0\rangle=1=\langle U_0,\hat\xi(e)\rangle.
\end{eqnarray*}
Then the associated cross ratio  is 
$$
b(f,v,e_0,u_0)=\langle\hat\xi^*(v),\hat\xi(f)\rangle
$$
From this expression it follows that $\chi^{n+1}=0$ and that by hyperconvexity $\chi^n\not=0$.
\subsubsection{Curves associated to cross ratios}
We now suppose that we have a strict cross ratio $b$ such that $\chi_b^{n+1}=0$. Assume also that for every pair  $e=(e_0,e_1,\ldots,e_n)$, $u=(u_0,u_1,\ldots,u_n)$
 of $n$-uples satisfying
\begin{eqnarray*}
j>i>0 \implies e_j\not= e_i\not= u_0,\ \  u_j\not=u_i\not=e_0,
\end{eqnarray*}
we have 
$$\chi_b^p(e,u)\not=0.$$
We  denote by $[\ldots : \ldots ]$ the projective coordinates on $\mathbb P(\mathbb R^n)$. Let $(e,u)$ be  as above.
We define two maps 
\begin{eqnarray*}
&\xi&\mapping{S^1\setminus\{u_0,e_0\}}{\mathbb P(\mathbb R^n)}{f}{(b(f,u_1,e_0,u_0):\ldots:b(f,u_n,e_0,u_0))}\\
&\xi^*&\mapping{S^1\setminus\{u_0,e_0\}}{\mathbb P(\mathbb R^n)}{v}{(b(e_1,v,e_0,u_0):\ldots:b(e_n,v,e_0,u_0))}
\end{eqnarray*}
We now prove that $b=b_{\xi,\xi^*}$. By Lemma \ref{unixi}, this will show that the curves $\xi$ and $\xi^*$ are unique up to projective transformations, and in particular do not depend on the choice of $e$ and $u$.

By our construction, 
\begin{eqnarray*}
\chi_{b_{\xi\xi^*}}^n(e,u)&=&\chi_{b}^n(e,u)\label{bbxi0}.
\end{eqnarray*}
Moreover, for every $f$ and $v$
\begin{eqnarray}
b(f,u_i,e_0,u_0)&=&b_{\xi,\xi^*}(f,u_i,e_0,u_0),\label{bbxi1}\cr b(e_i,v,e_0,u_0)&=&b_{\xi,\xi^*}(e_i,v,e_0,u_0)\label{bbxi2}.
\end{eqnarray}

We now choose  $e_n=f$, $u_n=v$. Developping the determinant $\chi^{n+1}(e,u)$ along the last line, the equation $\chi^{n+1}(e,u)=0$ leads to
$$
b(f,v,e_0,v_0)\chi_b^n(e,u)=F(\ldots,b(f,u_i,e_0,u_0),\ldots,b(e_i,v,e_0,u_0),\ldots)
$$
where the right hand term is a polynomial  in $b(f,u_i,e_0,u_0)$ and $b(e_i,v_i,e_0,u_0)$.
For the same reason we have
$$
b_{\xi,\xi^*}(f,v,e_0,v_0)\chi_b^n(e,u)=F(\ldots b_{\xi,\xi^*}(f,u_i,e_0,u_0),\ldots,b_{\xi,\xi^*}(e_i,v,e_0,u_0),\ldots).
$$
Therefore, using Equations (\ref{bbxi0}), (\ref{bbxi1}) and (\ref{bbxi2}) we obtain
\begin{eqnarray*}
b(f,v,e_0,v_0)\chi_b^n(e,u)&=&F(\ldots,b_{\xi\xi^*}(f,u_i,e_0,u_0),\ldots,b_{\xi\xi^*}(e_i,v,e_0,u_0),\ldots),\\
&=&b_{\xi\xi^*}(f,v,e_0,v_0)\chi_{b_{\xi\xi^*}}^n(e,u),\\&=&b_{\xi\xi^*}(f,v,e_0,v_0)\chi_{b}^n(e,u).
\end{eqnarray*}
It follows that $b(f,v,e_0,u_0)=b_{\xi,\xi^*}(f,v,e_0,u_0)$. Applying Relation (\ref{bir11}) twice, we easily obtain that $b= b_{\xi,\xi^*}$. Finally since $\xi^n(e,u)\not=0$, $\xi$ is hyperconvex as well as $\xi^*$. 

\section{The jet space $J^1(S^1,\mathbb R)$}\label{jetspacesect}
In this preliminary section, we aim to explain the geometry of  $J=J^1(S^1,\mathbb R)$,the space of one-jet of functions on $S^1$, and how the group $\cdiff$ acts on $J$. We then explain that this geometric description  characterises $J$. Finally, we  describe a homomorphism of $PSL(2,\mathbb R)$ in $\cdiff$, whose image acts faithfully and transitively on $J$.
\subsection{Description of the jet space}\label{jetspace}
We describe in this section   the geometric features of $J$ that we shall use later on, namely: 
\begin{itemize}
\item projections on $S^1$, $T^*S^1$ and $S^1\times\mathbb R$,
\item the action of $\cdiff$ on $J$,
\item a foliation of $J$ by affine leaves,
\item a flow,
\item a contact form.
\end{itemize}
\subsubsection{Projections}
Let $J=J^1(S^1,\mathbb R)$ be the space of one-jet of functions on $S^1$. We denote by $j^1_x(f)$ the one-jet of the function $f$ at the point $x$. We define the projections
\begin{eqnarray*}
\pi&:&\mapping{J}{S^1}{j^1_x(f)}{x}\cr
\pi_1&:&\mapping{J}{S^1\times\mathbb R}{j^1_x(f)}{(x,f(x))}\cr
\pi_2&:& \mapping{J}{T^*S^1}{j^1_x(f)}{d_xf}
\end{eqnarray*} 
We observe that each fibre of $\pi$ carries an affine structure.
\subsubsection{Action of $\cdiff$}\label{actionj1}
We denote by $Di\!f\!f^{h}(S^1)$ the group of $C^1$-diffeomorphisms of $S^1$ with Hölder derivatives. Let $C^{1,h}(S^1)$ the space of $C^1$-functions  on $S^1$ with Hölder derivatives.
The group $\cdiff$ acts naturally on $J$  by Hölder homeomorphisms in the following fashion: for every $C^1$-diffeomorphism $\phi$ of $S^1$, for every function $h$ in $C^1(S^1)$, we define the homeomorphism of $J$
$$
F=(h,\phi)\ : \ j^1_\theta(f)\mapsto j^1_{\phi(\theta)}((h+f)\circ\phi^{-1}).
$$ 
Alternatively, if we choose a coordinate $\theta$ on $S^1$ and consider the identification
$$
\mapping{J}{S^1\times\mathbb R\times\mathbb R}{j^1_\theta(f)}{(\theta, r,f)=(\theta,\frac{\partial f}{\partial \theta},f(\theta)),}
$$
then
\begin{eqnarray}
F=(h,\phi)\ : \ (\eta,r,f)\mapsto(\phi(\eta), \frac{\partial\phi}{\partial \theta}^{-1}\! \!\!\!(\eta)\big(\frac{\partial h}{\partial\theta}(\eta)+ r \big),f+h(\eta)).\label{fhfi}
\end{eqnarray}

The homeomorphism $F$ has the following properties
\begin{itemize}
\item $F$ preserves the fibres of $\pi$:$\pi\circ F=\phi\circ\pi$.
\item For every $x$ in $S^1$, $F$ restricted to $\pi^{-1}(x)$ is affine, in particular $C^\infty$, and the derivatives of $F\vert_{\pi^{-1}(x)}$ vary continuously on $J$.
\end{itemize}

\subsubsection{Foliation}
Let ${\mathcal F}$ be the 1-dimensional  foliation of $J$ given by the fibres of $\pi_1$.
We observe that the fibre through $z$ is identified with $T_x^*S^1$ for $x=\pi(z)$. This fibre also carries an affine structure invariant by the action of $\cdiff$ described above. 
\subsubsection{Canonical flow}
We define the {\em canonical flow} of $J$ to be the flow
$$
\phi_t(j^1(f))=j^1(f+t),
$$
where we identify the real number $t$ with the constant function that takes value $t$.
The canonical flow commutes with the action of $\cdiff$ on $J$. Notice also that 
$$
J/\phi_t=T^*S^1,
$$
and this identification turns $\pi_2$ in a principal $\mathbb R$-bundle.
\subsubsection{Contact form}\label{contactform}
We finally recall  that $J$ admits a contact form $\beta$.  If we choose a coordinate $\theta$ on $S^1$ and consider the identification
$$
\mapping{J}{S^1\times\mathbb R\times\mathbb R}{j^1_\theta(f)}{(\theta, r,f)=(\theta,\frac{\partial f}{\partial \theta},f(\theta))}
$$
then
$$
\beta= df - r d\theta.
$$  
Note that a legendrian curve for $\beta$ which is locally a graph above $S^1$ is the graph  of one-jet of a function. Moreover, the canonical flow $\phi_{t}$ preserves $\beta$.
Finally, we observe that $\beta$ is a connection form for the $\mathbb R$-principal bundle
 defined by $\pi_2$, whose curvature form is the canonical symplectic form of 
 $T^* S^1$. We prove.
 
\begin{proposition}\label{caractcontact}
Let $\alpha$ be a connection one-form on $J$. Assume that the curvature of $\alpha$ is the canonical symplectic form on $T^*S^1$.  Then there exists a diffeomorphism $\xi$ of $J$ which
\begin{itemize}
\item commutes with the action of the canonical flow,
\item preserves the fibres of $\pi$ and $\pi_2$,
\item is above a symplectic diffeomorphism of $T^*S^1$,
\item satisfies $\xi^*\beta=\alpha$.
\end{itemize}
\end{proposition}
\proof
We choose a coordinate $\theta$ of $S^1$. Thus, $T^*S^1$ is identified with $S^1\times \mathbb R$.  
Since $\alpha$ is a connection form there exists functions $\alpha_r$ and $\alpha_\theta$ such that
$$
\alpha=df +\alpha_r(r,\theta)dr + \alpha_\theta(r,\theta)d\theta.
$$
Let
$$
\gamma=\alpha-\beta=\alpha_r(r,\theta)dr + (\alpha_\theta(r,\theta)-r)d\theta.
$$
Since the curvature of $\alpha$ is $d\theta\wedge dr$,  $\gamma$ is closed. Therefore, there exist a function $h$ and a constant $\lambda$ such that 
$$
\gamma=dh+\lambda d\theta.
$$
Let $\mu_\lambda$ be the symplectic diffeomeorphism of $T^*S^1$ such that
$$
\mu_\lambda(\theta,r,f)=(\theta, r-\lambda, f).
$$
It is a straightforward check that
$$
\xi=\mu_\lambda\circ (h,id),
$$
satisfies the condition of our proposition. \qed
 
\subsection{A Geometric Characterisation of $C^\infty(S^1)\rtimes Di\!f\!f^\infty (S^1)$}
The following proposition suggests that $\mathcal F$, $\beta$ and $\phi_t$ characterise the action of $\cdiff$, or at least the smooth diffeomorphisms in $\cdiff$. 
\begin{proposition}\label{cidiff}
Let $\psi$ be a $\ci$-diffeomorphism of $J$. Assume that
\begin{enumerate} 
\item $\psi$ commutes with $\phi_t$, 
\item $\psi$ preserves $\beta$,
\item $\psi$ preserves $\mathcal F$,
\end{enumerate} 
then $\psi=(f,\phi)\in\cdiff$, where $f$ and $\phi$ are $\ci$.
\end{proposition}

In Corollary \ref{caractcdiff} we characterise
$\cdiff$.

\proof We us the notations and assumptions of the Proposition.
By (1) and (3), we obtain that $\psi$ preserves the fibres of $\pi$. Thus, there  exists a $\ci$-diffeomorphism $\phi$ of $T^*S^1$ such that
$$
\pi\circ\psi=\phi\circ\pi.
$$
Replacing $\psi$ by $\psi\circ(0,\phi^{-1})$, we may as well assume that $\phi=id$.
By (1), $\psi$ also preserves the fibres of $\pi_2$.

We choose coordinate in $S^1$ and from now on, we use the identification $J=S^1\times \mathbb R\times\mathbb R$ given in Paragraph \ref{contactform}. The previous
discussion shows that
$$
\psi (\theta,r,f)=(\theta, F(\theta,r), H(\theta,r,f)).
$$ 
Since $\psi^*\beta=\beta$, we obtain that
$$
dH - Fd\theta = df - rd\theta.
$$
Hence
$$
\frac{\partial H}{\partial f}=1,\ F-\frac{\partial H}{\partial \theta}=r, \ \frac{\partial H}{\partial r}=0.
$$
Therefore, there exists a function $g$ such that $H=f+g(\theta)$.
It follows that 
$$
\psi (\theta,r,f)=(\theta, r+\frac{\partial g}{\partial \theta}, f + g(\theta)).
$$
This exactly means that
$$
\psi(\theta,r,f)=(g,id). j^1(f).
$$
\qed

\subsection{$PSL(2,\mathbb R)$ and $C^\infty(S^1)\rtimes Di\!f\!f^\infty (S^1)$}

We consider the 3-manifold $\overline{J}=PSL(2,\mathbb R)$. It has a biinvariant Lorentz metric $g$ coming from the Killing form. We introduce the following objects.
\begin{itemize}
\item Let $\overline{\phi_t}$ be the one-parameter group of diagonal matrices acting on the right on $M$. 
\item Let $X$ be the vector field  generating $\overline\phi_{t}$, 
\item Let $\overline{\mathcal F}$ the  foliation by the right action of the one-parameter group of upper  triangular matrices $B$, 
\item Let $\overline\beta=i_X g$.
\end{itemize}
Then
\begin{proposition}\label{stan}
The one-form $\overline{\beta}$ is a contact form. Furthermore, there exists a $C^\infty$-diffeomorphism $\Psi$ from $J$ to $\overline{J}$, that sends $(\phi_t,{\mathcal F},\beta)$ to $(\overline{\phi_t}, \overline{\mathcal F},\overline\beta)$ respectively. Moreover, the diffeomorphism  $\Psi$ is unique up to left composition by a $\ci$-diffeormorphism element of $\cdiff$.
\end{proposition}
As an immediate application, we have
\begin{definition}\label{standef}
By  Propositions \ref{cidiff} and \ref{stan} and since  $PSL(2,\mathbb R)$ preserves $\overline{\phi_t}$, $\overline{\mathcal F}$ and 
$\overline\beta$,  we have a group homomorphism 
$$
\iota:\ \mapping{PSL(2,\mathbb R)}{C^\infty(S^1)\rtimes Di\!f\!f^\infty(S^1)}{g}{\Psi\circ g\circ \Psi^{-1},}$$ well defined up to conjugation by  a $\ci$-diffeormorphism element of $\cdiff$. The corresponding representation is called {\em standard}.
\end{definition}

\vskip 0.5truecm
\proof  It is immediate to check that
$$
PSL(2,\mathbb R)\rightarrow PSL(2,\mathbb R)/\overline\phi_t=W,
$$
is a $\mathbb R$ principal bundle, whose connexion form is $\overline{\beta}$ and whose curvature is 	 symplectic.  It follows that $\overline\beta$ is a contact form.

Let $\overline{\pi}$ be  the projection 
$$
W\rightarrow {\mathbb{RP}}^1=PSL(2,\mathbb R)/B.
$$
We observe that $W$ is diffeomorphic to the annulus $S^1\times\mathbb R=T^*S^1$ in such a way $\overline{\pi}$ coincides with the projection $\pi$ on the first  factor. We may actually choose this diffeomorphism $\psi_0$ to be symplectic, still sending the fibres of $\overline{\pi}$ to the fibres of $\pi$. 

Let now $\xi$ be the  symplectic diffeomorphism of $T^*S^1$ obtained by Proposition \ref{caractcontact} applied to $\alpha=(\psi_0^{-1})^*\overline\beta$. It follows $\Psi=\xi\circ\psi_0$ has all the properties required. Finally, by Proposition \ref{cidiff}, $\Psi$ is well defined up to multiplication with an element of 
$C^\infty(S^1)\rtimes Di\!f\!f^\infty (S^1)$.

Since the right action of $PSL(2,\mathbb R)$ preserves $\overline{\mathcal F}$, $\overline{\beta}$, and $\overline{\phi}_t$, we obtain a representation of $PSL(2,\mathbb R)$ well defined up to conjugation
$$
\mapping{PSL(2,\mathbb R)}{C^\infty(S^1)\rtimes Di\!f\!f^\infty (S^1)}{g}{\Psi\circ g\circ\Psi^{-1}}
$$
\qed
\vskip 1 truecm
{\sc Remark :} The action of $PSL(2,\mathbb R)$  on ${\mathbb {RP}}^1$ gives rise to an  embedding of $PSL(2,\mathbb R)$ in $Di\!f\!f^\infty (S^1)$. The above standard representation is a nontrivial  extension of this representation. Indeed, the natural lift of the action of  $Di\!f\!f^\infty (S^1)$ on $T^*S^1$ does not act transitively  since it preserves the zero section. On the contrary, $PSL(2,\mathbb R)$ does act transitively through the standard representation.

\section{Homomorphisms of $\grf$ in $\cdiff$}\label{homosect}

Our aim in this section is to describe a "good" set of homomorphisms of $\grf$ in $\cdiff$. For a semi-simple Lie group $G$, a "good" set of homomorphisms of $\grf$ in $G$ is the set  of reductive homomorphisms ({\it i.e.} such that the Zariski closure of the image is reductive) it satisfies the following properties :
\begin{itemize}
\item it is open in the set of all homomorphisms,
\item $G$ acts properly on it.
\end{itemize}.
As a consequence, the set of "good" representations, that is the set of "good" homomorphisms up to conjugacy, is a Haussdorff space. 

We define "good" (actually {\it $\infty$-Hitchin}) homomorphisms in Paragraph \ref{goodhomo} and show in Theorem \ref{maingoodhomo}, that they satisfy analogous properties to those of reductive homomorphisms in semi-simple Lie groups that we just described. We also show that to a  $\infty$-Hitchin representation is associated a cross ratio. In Theorem \ref{act}, we finally explain how two representations with the same cross ratio are related.

\subsubsection{Fuchsian and $H$-Fuchsian homomorphims}

\begin{definition}{\sc [Fuchsian homomorphism]}
Let $\rho$ be a Fuchsian homomorphism of $\grf$ in $PSL(2,\mathbb R)$.  We say that the composition of $\rho$ by the standard representation $\iota$ ({\em cf} Definition \ref{standef}) is {\em  $\infty$-Fuchsian}, or in short {\em Fuchsian} when there is no ambiguity.
\end{definition}

\begin{definition}{\sc [$H$-Fuchsian action on $S^1$]}\label{defhomh}
We  say an action of $\grf$ on $S^1$  is {\em $H$-Fuchsian} if it is an action by $C^1$-diffeomorphisms with Hölder derivatives, Hölder conjugate to the action of a cocompact group of $PSL(2,\mathbb R)$ on $\pr$.
\end{definition}
\subsubsection{Anosov and $\infty$-Hitchin representations}\label{goodhomo}
We begin with a general definition.
\begin{definition}
\label{defcontract} For a foliation $\mathcal F$ on a compact space, we denote by  $d_{\mathcal F}$ a Riemannian distance along the leaves of $\mathcal F$ which comes from a leafwise continuous metric. In particular, 
$$
d_{\mathcal F}(x,y)<\infty {\hbox{ iff $x$ and $y$ are in the same leaf}}.
$$
We say a  {\em flow $\phi_{t}$ contracts uniformly the leaves of a foliation $\mathcal F$}, if $\phi_t$ preseves $\mathcal F$ and if moreover 
$$
 \forall \epsilon >0, \ \exists \alpha >0,\  \exists t_{0}: \ t \geq t_{0}, \ d_{\mathcal F}(x,y)\leq\alpha \implies  d_{\mathcal F}(\phi_{t}(x),\phi_{t}(y))\leq\epsilon
$$
This definition does not depend on the choice of $d_{\mathcal F}$.
\end{definition}
We now use the notations of Section \ref{jetspace}. 
\begin{definition}{\sc [Anosov homomorphims]}
A homomorphism $\rho$  of $\grf$ in $\cdiff$ is {\em Anosov} if:
\begin{itemize}
\item $\rho(\grf)$  acts with a compact quotient on $J=J^1(S^1,\mathbb R)$,
\item The flow induced by the canonical flow $\phi_t$ on $J/\rho(\grf)$ contracts uniformly the leaves of $\mathcal F$ ({\it cf.} definition above)
\item The induced action on $S^1$ is $H$-Fuchsian.
\end{itemize}
We denote the space of Anosov homomorphims by $\Hom^*$.
\end{definition}
\vskip 1 truecm
\noindent\rmks
\begin{enumerate}
\item If  $\rho$ is Anosov, one should observe that in general $\rho(\grf)$ only acts by homeomorphisms on $J$. However, if  we consider $J$ as a $C^\infty$-filtrated space  (See the definitions in the Appendix \ref{aff}), whose nested foliated structures are given by the fibres of the two projections $\pi_2$ and $\pi$, then  $\rho(\grf)$
 acts by $C^\infty$-laminated maps ({\it i.e} smoothly along the leaves with continuous derivatives). Therefore, $J/\rho(\grf)$ has the structure of a $C^\infty$-lamination such that the flow induced by $\phi_t$ is  $C^\infty$-leafwise.  
\item $\infty$-Fuchsian representations of $\grf$ in
$\cdiff$ provide our first examples of elements of Anosov homomorphisms.
\end{enumerate}

\begin{definition}{\sc [$\infty$-Hitchin homomorphims]}
Let  $\Hom_H$ be the connected component of $\Hom^*$ containing the Fuchsian homomorphisms. Elements of $\Hom_H$ are called{ \em $\infty$-Hitchin homomorphisms}.
\end{definition}

The group $\cdiff$ acts by conjugation on $\Hom_H$. The canonical flow $\phi_t$ is in the center of $\cdiff$ and hence acts trivially on $\Hom_H$. It follows that $\cdiff/\phi_t$ acts on $\Hom_H$.

The main result of this section is the following: 

\begin{theorem}\label{maingoodhomo}
The set 
$\Hom_H$ is open in the space of all homomorphisms. Furthermore, the action of $\cdiff/\phi_t$ by conjugation on  $\Hom_H$  is proper and the quotient is Haussdorff
\end{theorem}

We prove later on a refinement (Theorem \ref{act}) of the second part of this result.
We shall denote by $\Rep_H$ the quotient $\Hom_H/\cdiff$ where the right action is by conjugation.

The proof of the Theorem falls into two parts.  We first show that  $\Hom_H$ is open. Then, we prove that $\cdiff$ acts properly  with a Haussdorff quotient. In this proofs, we  use the definitions and propositions obtained in the Appendix A (\ref{applam}).
Accordingly, from now on, we shall consider $J$ as a $C^\infty$-filtrated space, whose nested foliated structures are given by the fibres of the two projections $\pi_2$ and $\pi$.

\subsection{Openness of $\Hom^*$ and a Stability Lemma}

Let $\rho$ be a representation in $\cdiff$. We denote by $\overline\rho$ the associated representation in $Di\!f\!f^{h}(S^1)$.

The following Stability Lemma  implies immediately the openness of $\Hom^*$. Moreover, using corollaries of this Lemma, we  associate to every $\infty$-Hitchin representation a cross ratio and a spectrum in Paragraph \ref{crosshomh} and \ref{specthomh}.

\begin{lemma}\label{conjug0}
Let $\rho_0$ be an Anosov representation. Then for $\rho$ close enough to $\rho_0$,  there exists a Hölder homeomorphism $\Phi$ of $J$ close to the identity which is a $C^\infty$-filtrated immersion as well as its inverse such that
$$
\forall \gamma \in \grf, \rho_0(\gamma)=\Phi^{-1}\circ\rho(\gamma)\circ\Phi.
$$
\end{lemma}

\subsubsection{Minimal action on the circle}\label{minactsect}

We prove the following preliminary Lemma which is independent of the rest of the article.

\begin{lemma}\label{minimal}
Let $\Gamma$ be a group acting on $S^1$ by homeomorphisms. Suppose that every orbit of the action on $S^1\times S^1$ is dense. Suppose that every element of $\Gamma$ different from the identity has exactly two fixed points, one attractive the other repulsive. Let $\rho$ be a representation of  $\Gamma$ in the group $\rm{Homeo}(S^1)$ of homeomorphisms of $S^1$. Let $f$ be a continuous map of degree different from zero, from $S^1$ to $S^1$ such that
$$
\forall \gamma\in\Gamma,  f\circ \gamma=\rho(\gamma)\circ f.
$$ 
Then $f$ is a homeomorphism.
\end{lemma}

\proof Since $f$ has nonzero degree, it is surjective.

We now prove $f$ is injective. Let $a$ and $b$ two points of $S^1$ such that $a\not=b$ and $f(a)=f(b)=c$. Let $I$ and $J$ be the connected components of $S^1\setminus\{a,b\}$.

Since $f$ has a nonzero degree, either $S^1\setminus \{c\}\subset f(I)$ 
or $S^1\setminus \{c\}\subset f(I)$. Assume  $S^1\setminus \{c\}\subset f(I)$.  By the density of orbits on $S^1\times S^1$, we can find an element $\gamma$ in $\Gamma$ such that $\gamma^+\in I$ and $\gamma^-\in J$.

We observe that for all $n$, $f(\gamma^n(I))=S^1\setminus \{\rho(\gamma^n)(c)\}$.
It follows we can find $\epsilon >0$, and two sequences $\{x_n\}_{n\in\mathbb N}$, $\{x_n\}_{n\in\mathbb N}$ such that
\begin{itemize}
\item $x_n\in \gamma^n(I)$, $y_n\in\gamma^n(I)$,
\item $d(f(x_n),f(y_n))>\epsilon$.
\end{itemize}
The contradiction follows since 
$$
\lim_{n\rightarrow\infty}d(x_n,y_n)=0.
$$
\qed

\subsubsection{Proof of Lemma \ref{conjug0}}
\proof
Let   $\rho_0$ be an an element of $\Hom_H$. 
$$
\rho_0:\grf\rightarrow\cdiff.
$$
By assumption, $P=\rho_0(\grf)\backslash J$ is compact. The topological space $P$ is a $C^0$-manifold which is $C^\infty$-laminated by the fibres of the projection $\pi: J\rightarrow S^1$ and $\pi_2: J\rightarrow T^*S^1$, since $\cdiff$ acts as $C^\infty$-laminated maps on $J$. The covering $J\rightarrow P$ is a Galois covering.

Let now $\rho$ be a representation of  $\grf$ close enough to $\rho_0$. Let $\Gamma=\rho_0(\grf)$. According to Theorem \ref{laminholo},  we obtain a $\rho\circ\rho_0^{-1}$-equivariant $C^\infty$-laminated immersion $\Phi$ from 
$J$ to $J$, arbitrarily close to the identity on compact sets. It remains to show that $\Phi$ is a homeomorphism.

Since $\Phi$ is a laminated map, we obtain a Hölder map $f$ close to the identity from $S^1$ to $S^1$ such that
$$
\pi\circ\Phi=f\circ\pi,
$$
and
$$
\forall \gamma\in\Gamma, f\circ\overline{\rho}_0(\gamma))=\overline{\rho}(\gamma) \circ f
$$
We recall that $\rho_0$ is Anosov, therefore $\overline\rho_0$ is a $H$-Fuchsian representation, and satisfies the hypothesis of Lemma \ref{minimal}. By this Lemma, and since $f$ being close to the identity is not of degree 0, we obtain that  $f$ is a homeomorphism. Thus,  $\overline{\rho_1}$ is also $H$-Fuchsian.
 
This laminated equivariant immersion $\Phi$ induces an affine structure on the leaves of of $P$, according to the definition of Paragraph \ref{aff}. By compactness of $P$, we deduce there exists an action of a one-parameter group $\psi_t$ such that, on $J$
$$
\Phi\circ\psi_t=\phi_t\circ\Phi,
$$
where $\phi_t$ is  the canonical flow. We observe that $\psi_t$ preserves any given leaf and acts as a one-parameter group of translation on it.

Recall that the canonical flow contracts uniformly the leaves of $\mathcal F$ for on $J/\rho_0(\grf)$. It follows  the same holds for $\psi_t$, for $\Phi$ close enough to the identity.

By Lemma (\ref{affcomp}), we obtain that the affine structure  induced by $\phi$ is leafwise complete. In particular, for every $x$ in $S^1$, $\Phi$ is an affine bijection from $\pi^{-1}\{x\}$ to  $\pi^{-1}\{f(x)\}$. Since $f$ is a homeomorphism, we deduce that $\Phi$ itself is a homeomorphism, its inverse being also a filtrated immersion.

This concludes the proof. \qed 
\vskip 0.5truecm

We have the following immediate Corollary, 

\begin{coro}\label{conjug1} Let $\rho_0$ and $\rho_1$ be two elements of $\Hom_H$. Then there exists a Hölder homeomorphism $\phi$ of $J$ which is a filtrated immersion as well as its inverse such that
$$
\forall \gamma \in \grf, \rho_0(\gamma)=\phi^{-1}\circ\rho_1(\gamma)\circ\phi.
$$
\end{coro}
Recall that $\cdiff$ acts on $T^*S^1=J/\phi_t$. We finally prove 
\begin{proposition}\label{conjug1cross}
Let $\rho$ be an element in $\Hom_H$. Then there exists a unique $\rho(\grf)$-equivariant Hölder homeomorphism  $\Theta_\rho$  from $T^*S^1$ to
$$\bgrf^{2*}=\bgrf^2\setminus\{(x,x)/x\in\bgrf\}$$ such that
$$ \exists f : S^1\rightarrow\bgrf,  \Theta_\rho(T^*_x S^1)\subset \{f(x)\}\times\bgrf.$$
\end{proposition}
\proof By Corollary \ref{conjug1}, it suffices to show this result for the action of a cocompact subgroup of $PSL(2,\mathbb R)$. In this  case the statement is obvious. \qed

\subsubsection{Cross ratio associated to elements of $\Hom_H$}\label{crosshomh}
We show that we can asociate a cross ratio to every element $\rho$ of $\Hom_H$.

Let $(a,b,c)$ be  a triple of distinct elements of $\bgrf$. We define $[b,c]_a$ to be the closure of the connected component of $\bgrf\setminus\{b,c\}$ not containing $a$. 
Let $q=(x,y,z,t)$ be a quadruple of elements of $\bgrf$, we define $\widehat q$ to be following closed curve, embedded in $\bgrf^{2*}$,
$$
\widehat q=(\{x\}\cup [y,t]_x)\cup([x,z]_t\cup\{t\} )\cup(\{z\}\cup [y,t]_z)\cup([x,z]_y\cup\{y\} ).$$
We choose the orientation on $\widehat q$ such that $(x,y),(x,t),(z,t),(z,y)$ are cyclically oriented. We define $\epsilon_q\in\{-1,1\}$ to be the sign of the classical cross ratio of $q$ on $\bgrf$.

Let $\Theta_\rho$ be the map from  $T^*S^1$ to
$\bgrf^{2*}$ as defined in Proposition \ref{conjug1cross}. 
\begin{definition}{\sc [Associated cross ratio]}
We define the {\em associated cross ratio} to $\rho$ by 
$$
\widehat b_\rho(x,y,z,t)={\epsilon_q}\exp\big(\frac{1}{2}{\int_{\Theta_\rho^{-1}(\widehat q)}r d\theta}\big).
$$
\end{definition}
\vskip 1 truecm
\noindent\rmk We observe that the cross ratio just depends of the action of $\grf$ on $T^*S^1$. Here is another formulation of this observation.  Let $\Omega^h(S^1)$ be the space of Hölder one-form on $S^1$. Note that $\ocdiff$ naturally acts on $T^*S^1$. We also have a natural homorphism 
$$d:\mapping{\cdiff}{\ocdiff}{(f,\phi)}{(df,\phi)},$$ whose kernel is the canonical flow.
Then two representations $\rho_1$ and $\rho_2$ such that $d\rho_1=d\rho_2$ have the same associated cross ratio. In other words, the cross ratio only depends on  the representation as with values in $\ocdiff$,
\subsubsection{Spectrum}\label{specthomh}
\begin{definition}{\sc [$\rho$-length]}
Let $\rho$ be a representation in $\Hom_H$. Let  
$\gamma$ be an element  in $\grf$,  the {\em $\rho$-length} of $\gamma$, $l_{\rho}(\gamma)$ is the positive number $t$ such that
$$
\exists u\in J, \phi_t(u)=\rho(\gamma) u.
$$
\end{definition}
We observe that the existence and uniqueness of such a number $t$ is follows by Corollary (\ref{conjug1}) and the description of standard representations.

In other words, $l_{\rho}(\gamma)$ is the length of the periodic orbit of $\phi_t$ in $J/\rho(\grf)$ freely homotopic to $\gamma$. It is clear that  $l_{\rho}(\gamma)$ just depends on the conjugacy class of $\gamma$.
\begin{definition}{\sc [Spectrum]}
The {\em marked spectrum} of $\rho$ is the map
 $$
 l_{\rho}:\gamma\mapsto l_{\rho}(\gamma).
 $$
We say a representation is {\em symmetric} if $l_\rho(\gamma)=l_\rho(\gamma^{-1})$. 
\end{definition}
\subsubsection{Compatible representations and Ghys deformations}

\begin{definition}{\sc [Compatible representation]}\label{compa} We say an element $\rho$ is {\em compatible} if its marked spectrum coincides with the periods of the associated cross ratio. 
\end{definition}
We observe that every compatible representation is symmetric.
Conversely, we prove in Proposition \ref{width} that if $l_{b_\rho}$ is the period of the cross ratio associated to the representation $\rho$, then
$$
l_{b_\rho}(\gamma)=\frac{1}{2}(l_\rho(\gamma)+l_\rho(\gamma^{-1})).
$$
Consequently, every symmetric representation is compatible.

We shall later obtain lot's of compatible representations. 
\begin{definition}{\sc [Ghys deformation]}
If $\rho$ is a compatible representation, if $\omega$ is a small enough nontrivial element in $H^1(\Gamma,\mathbb R)$, then 
$\rho^{\omega}$ defined by
$$
\rho^{\omega}(\gamma)=e^{\omega(\gamma)}.\rho(\gamma).
$$
is an element of $\Hom_H$. We call such  a representation $\rho^{\omega}$ a {\em Ghys deformation} of $\rho$.
\end{definition}
We observe that  $\rho^{\omega}$ and $\rho$ have the same associated cross ratio since they have the same action on $T^*S^1$. 
However
$$
l_{\rho^\omega}(\gamma)=l_\rho(\gamma)+\omega(\gamma).
$$
It follows that two representations can have the same cross ratio but different marked spectrum.

\subsection{Action of $\cdiff$ on $\Hom_{H}$}
We prove in this section that $\cdiff/\phi_t$ acts properly on $\Hom_H$ and more precisely 

\begin{theorem}\label{act}
The group $\cdiff/\phi_t$ acts  properly on $\Hom_H$ with a Haussdorf quotient.  

Moreover, if 
two representations $\rho_0$ and $\rho_1$ elements of  $\Hom_H$  have the  same associated cross ratio, then there exists $\omega\in H^1(\grf)$ such that $\rho_0$ and $\rho_1^\omega$ are conjugated by some element in $\cdiff$.

Consequently two representations with the same cross ratio and the same spectrum are conjugated in $\cdiff$.
\end{theorem}
\noindent\rmk We shall see later that two representations with the same spectrum are not necessarily conjugated

\subsubsection{A preliminary lemma}
Let $\Omega^h(S^1)$ be the space of Hölder one-form on $S^1$. We observe first that $\ocdiff$ acts naturally on $T^*S^1$ and preserves the area.
\begin{lemma}\label{areapreserve}
Let $G$ be an area preserving homeomorphism of $T^*S^1$. Assume there exists a homeomorphism  $f$ such that 
$$
G(T_x S^1)\subset T_{f(x)}S^1.
$$
Then $G$ belongs to $\ocdiff$.
\end{lemma}
\proof We use the coordinates $(r,\theta)$ on $T^* S^1$. By hypothesis,
$$
G(r,\theta)=(g(r,\theta),f(\theta)).
$$
Since $G$ preserves the area, we obtain
\begin{eqnarray}
(r_0-r_1)(\theta_0-\theta_1 ) &=&\int_{f(\theta_0)}^{f(\theta_1)}(\int_{g(r_0,\theta)}^{g(r_1,\theta)}dr)d\theta \\
&=&\int_{f(\theta_0)}^{f(\theta_1)}(g(r_0,\theta)-g(r_1,\theta))d\theta \label{eqocdiff}\end{eqnarray}
Hence  $g(r,\theta)$ is affine in $r$:
$$
g(r,\theta)= \omega(\theta) +r\beta(\theta),
$$
where $\omega(\theta)$ and $\beta$ are Hölder.
Since $G$ is a homeomorphism, we observe that for all $\theta$, $\beta(\theta)\not=0$.
By Equation \ref{eqocdiff}, we obtain that
$$
f^{-1}(\theta_0)-f^{-1}(\theta)=\int_{\theta_0}^{\theta_1}\beta(\theta)d\theta.
$$
It follows that $f^{-1}$ is in $C^{1,h}(S^1)$ with $df^{-1}=\beta$. Since $\beta$ never vanishes, $f$ is actually is a diffeomorphism, and $G$ belongs to $\ocdiff$. \qed

\subsubsection{Conjugation in $\ocdiff$}
Let $\rho_0$ and $\rho_1$ be two representations of $\grf$ in $\cdiff$ with same associated cross ratio.

Let $\Omega^h(S^1)$ be the space of Hölder one-form on $S^1$. We recall that $\ocdiff$ acts naturally on $T^*S^1$ and that $\cdiff$ naturally maps to $\ocdiff$. We first prove

\begin{proposition}\label{conjugocdiff}
Let $\rho_0$ and $\rho_1$ be two representations of $\grf$ in $\ocdiff$ with same associated cross ratio.  Then there exists a unique  element $G=G^{\rho_0,\rho_1}$ in $\ocdiff$ such that
$$
\forall \gamma\in \grf, \forall y\in T^*S^1,\ \ G(\rho_0(\gamma)(y))=\rho_1(\gamma)(G(y)).
$$

Moreover, for  any representation $\rho$,  for any  neighbourhood $V$ of the identity in $\ocdiff$, there exists a neighbourhood $U$ of $\rho$, such that  if two representations $\rho_1$ and $\rho_0$  in $U$ have the same cross ratio,   then
$$
G^{\rho_0,\rho_1}\in V.
$$
\end{proposition}
\proof Applying Proposition \ref{conjug1cross} twice, we obtain a unique Hölder homeomorphism $G$  of $T^*S^1$ and a homeomorphim $f$ of $S^1$ such that
\begin{eqnarray*}
\forall \gamma\in \grf, \ \ G\circ \rho_0(\gamma)&=&\rho_1(\gamma)\circ G,\\
G(T_x S^1)&\subset& T_{f(x)}S^1.
\end{eqnarray*}
If $\rho_0$ and $\rho_1$ have the same associated cross ratio, then $G$ preserves the area. By Lemma \ref{areapreserve}, $G$ belongs to $\ocdiff$. This concludes the proof of the first part of the Proposition. The second part follows from Lemma \ref{conjug0}.

\subsubsection{Second step : from $\ocdiff$ to $\cdiff$}\label{proofact5}
We prove the following Lemma
\begin{lemma}\label{act5}
Let $\phi:\gamma\rightarrow\phi_\gamma$ be a  representation of $\grf$ in $Di\!f\!f^h (S^1)$ with nonzero Euler class. Let $\alpha$ a a continuous one-form. Let 
$$
\mapping{\grf}{C^1(S^1)}{\gamma}{f_\gamma}
$$
Assume that
\begin{eqnarray*}
\phi_\gamma^*(f_\eta)+f_\gamma&=&f_{\eta\gamma},\\
\phi_\gamma^*(\alpha)-\alpha&=&df_\gamma.
\end{eqnarray*}
Then $\alpha$ is exact.
\end{lemma}
\proof Let 
$$
\kappa=\int_{S^1}\alpha.
$$
We write 
$S^1=\mathbb R/\mathbb Z$. We lift all our data to $\mathbb R$, and use $\tilde u$ to describe a lift  of $u$. In particular, we have
$$
\tilde\phi_{\gamma\eta}(x)=\tilde\phi_\gamma\circ\tilde\phi_\eta (x)+c(\gamma,\eta)\in\mathbb Z,
$$
where
$$
c\in H^2(\grf,\mathbb Z),
$$
is a representative of the Euler class of the representation $\phi$.
We now write $\tilde \alpha=d h$. We observe that
$$
\forall m\in\mathbb Z, h(x+m)=h(x)+m\kappa.
$$
Let
$$
\eta :\grf\rightarrow {\mathbb R},
$$
be the function such that
$$
\tilde\phi_\gamma^*(h)-h= \tilde f_\gamma+ \eta (\gamma).
$$
We obtain
\begin{eqnarray*}
c(\gamma,\eta)\kappa&=&h\circ\tilde\phi_{\gamma\eta}-h\circ\tilde\phi_\gamma\circ\phi_\eta\\
&=&\tilde\phi_{\gamma\eta}^*(h) -\tilde\phi_\eta^*\tilde\phi_\gamma^*(h)\\
&=&(\tilde\phi_{\gamma\eta}^*(h)-h) -(\tilde\phi_\eta^*\tilde\phi_\gamma^*(h)-\tilde\phi_\eta^*(h))-(\tilde\phi_\eta^*(h)-h)\\
&=&\tilde f_{\gamma\eta}-\tilde\phi_\eta^*\tilde f_\gamma -\tilde f_\eta\\
&+& \eta (\gamma\eta)-\eta (\gamma) -\eta (\eta),\\
&=& \eta (\gamma\eta)-\eta (\gamma) -\eta (\eta).
\end{eqnarray*}
But, since the Euler class is nonzero in cohomology by hypothesis, we obtain that $\kappa=0$, hence $\alpha$ is exact.\qed

\subsubsection{Proof of Theorem \ref{act}}

We first prove that if two representations have the same cross ratio then, after a Ghys deformation, they are conjugate by an element of $\cdiff$.

Let $\rho_0$ and $\rho_1$ be two elements of $\Hom_H$ with the same cross ratio. Write   
$$
\rho_i(\gamma)=(f^i_\gamma,\phi^i(\gamma)).
$$
By Proposition \ref{conjugocdiff}, we know that, as representations in $\ocdiff$ they are conjugated in $\ocdiff$, by some element $G=(\alpha,F)$, where $F\in Di\!f\!f^{h}(S^1)$.

In particular, after conjugating by $(0,F)\in\cdiff$, we may assume that $$
\phi^0_\gamma=\phi^1_\gamma:=\phi_\gamma.
$$
Since $\rho_i$ is a representation
$$
\phi_\eta^*(f^i_\eta)+f_\gamma=f_{\eta\gamma}.
$$
By Proposition \ref{conjugocdiff}, there exists $\alpha\in \Omega^h(S^1)$
such that
$$
\phi_{\gamma}^*\alpha -\alpha= d f^1_\gamma -d f^0_\gamma.
$$
Applying Lemma \ref{act5}, to 
$$
f_\gamma=f^1_\gamma -f^0_\gamma,$$
we obtain that $\alpha=dg$.

Conjugating by $(g,Id)$, we may now as well assume that $\alpha=0$.
It follows that there exists
$$
\omega:\grf\rightarrow \mathbb R,
$$
such that
$$
f^1_\gamma -f^0_\gamma = \omega(\gamma).
$$
This exactly means that $\rho_0=\rho_1^\omega$.

This proves that  two representations with the same cross ratio are conjugated after a Ghys deformation. Consequently two representations with  the same cross ratio and spectrum are conjugated.

Finally, the proof of properness goes as follows.
Suppose that we have a sequence of representation $\{\rho_n\}\inn$ converging to $\rho_0$, a sequence of elements $\{\psi_n\}\inn$ in $\cdiff$ such that $\{\psi_n^{-1}\circ \rho_n\circ\psi_n\}\inn$ converges to $\rho_1$. Let $d$ be the homomorphism $\cdiff\rightarrow\ocdiff$. We aim to prove that $\{d\psi_n\}\inn$ converges.

Since  $\rho_n$ and $\psi_n^{-1}\circ \rho_n\circ\psi_n$ have the same spectrum and cross ratio, it follows that $\rho_1$ and $\rho_0$ have the same spectrum and cross ratio. From the first part of the proof it follows that there exists $F$ such that
$$
\rho_0=F^{-1}\circ \rho_1\circ F.
$$ 
We may therefore as well assume that $\rho_0=\rho_1$.

From Proposition \ref{conjugocdiff}, $\{d(\psi_n)\}\inn$ converges to the identity.\qed 

\section{A Conjugation Theorem}\label{conjugtheosect}

We will work in the following setting. 
\begin{enumerate}
\item
Let $\kappa$ be  a Hölder homeomorphism of $S^1$. We write
$$
D_{\kappa}=S^1\times S^1\setminus\{(s,t)/\kappa(s)=t\}.
$$
\item Let $p: M\rightarrow D_{\kappa}$ be a principal $\mathbb R$-bundle over $D_{\kappa}$ equipped
with a connection $\nabla$.  
\end{enumerate}
We define
\begin{itemize}
\item  Let $\mathcal L$ be the foliation of $M$ by parallel vector fields along the curves $c_s: t\mapsto(s,t)$.
\end{itemize}
We suppose
\begin{itemize}
\item Let $\omega$ be  the  curvature $\omega$ of $\nabla$. Let $f$ be such that  $\omega=f(s,t)ds\wedge dt$.  We suppose that $f$ is positive and Hölder. 
\end{itemize}
In this section, we aim to prove the following Theorem. 

\begin{theorem}\label{conjugtheo}
Assume that $\grf$ acts on $M$ by $C^1$-diffeomorphisms with Hölder derivatives. We suppose that this action preserves $\nabla$, and  that 
\begin{enumerate}
\item $M/\grf$ is compact. \label{cond1}
\item The action of $\mathbb R$ on $M/\grf$ contracts the leaves of $\mathcal L$ ({\em cf} definition \ref{defcontract}).\label{cond2}
\item There exists two $H$-Fuchsian  representations $\rho_{1}$ and $\rho_{2}$ of $\grf$, with $\rho_2=\kappa^{-1}\rho\kappa$ such that
$$
\forall\gamma\in\grf,\ \  p(\gamma u)= (\rho_{1}(\gamma)(p(u)),\rho_{2}(\gamma)(p(u))).
$$  
\end{enumerate}
Then there exists homeomorphism $\hat\psi$ from $M$ to $J$,   unique up to right composition by an element of $\cdiff$ an area preserving 
homeomorphism $\psi$ from $D$ to $T^*S^1$ such that
\begin{itemize}
\item
$
\pi_2\circ\hat\psi=\psi\circ p,
$ 
\item $\hat\psi$ commutes with the $\mathbb R$-action,
\item $\hat\psi$ sends $\mathcal L$ to $\mathcal F$,
\end{itemize}
and a  representation $\rho$ of $\grf$ in $\cdiff$ such that
$$ \hat\psi\circ\gamma=\rho(\gamma)\circ\hat\psi.$$
In particular, $\rho$ is Anosov. 

Finally, assume that our data $\kappa$, $\nabla$ depends continuously on a parameter. Then, we can choose $\hat\psi$ to depend continuously on this parameter.

\end{theorem}

In the first three sections we give a geometric description of elements of $\cdiff$ by their action on $J$. Then, we conclude the proof in the last  two sections.

To simplify the notations, we fix a trivialisation of $M=\mathbb R\times D_{\kappa}$ so that the vector fields  $U_{(s,u)} : t\mapsto (u,s,t)$ are parallel along $c_s$. 
Let 
\begin{eqnarray*}
\pi&:&\mapping{M}{S^1}{(u,s,t)}{t,}\cr
\pi_1&:&\mapping{M}{\mathbb R\times S^1}{(u,s,t)}{(u,s),}\cr
p &:& \mapping{M}{D_{\kappa}}{(u,s,t)}{(s,t).}
\end{eqnarray*} 
In particular, the foliation $\mathcal L$ is the foliation by the fibres of $\pi_1$.

For later use, using the above notations, we  state the following obvious result which allows us to restate Condition \ref{cond2} of Theorem \ref{conjugtheo}.

\begin{proposition}\label{conditionalternate}
Suppose $M/\pi_1(S)$ is compact. Then the following conditions (2) and (2') are equivalent 
\begin{itemize}
\item[(2)] The action of $\mathbb R$ on $M/\grf$ contracts the leaves of $\mathcal L$.
\item[(2')]
Let $\seq{t}$ be a sequence of real numbers going to $+\infty$. Let $\seq{\gamma}$ be a sequence of elements of $\grf$. Let $(u,s,t)$ be an element of $M$ such that 
$$
\{\gamma_m (u+t_m,s,t)\}_{m\in\mathbb M},
$$
converges to $(u_0,s_0,t_0)$ in $M$, then for all $w$ in $S^1$, with $w\not=\kappa(s)$, then
$$
\{\gamma_m (s,w)\}_{m\in\mathbb M},
$$
converges to $(s_0,t_0)$.
\end{itemize}
\end{proposition}

\subsection{$\pi$-exact symplectic Homeomorphisms of the Annulus}
We aim to describe in a geometric way the group  $\cdiff$ generalising Proposition \ref{stan}. Roughly speaking, we will describe it as  a subgroup of the central extension of exact symplectic "homeomorphisms" of the Annulus (i.e. $T^*S^1$). However, the notion of exact symplectic homeomorphism does not make sense in an obvious way, but some homeomorphisms may be coined as "exact symplectic" as we shall see.

\subsubsection{Line bundle}
We consider the real vector line bundle $p: L\mapsto \mathbb T^*S^1$ equipped with a connection $\nabla$ of curvature $\omega$. If $\beta$ is the Louville form, we identify $L$ with $T^*S^1\times \mathbb R$, with the connection
$$
\nabla=D +\beta.
$$
It is well known, and we shall recall the construction in Section \ref{sympham}, that the action of every Hamiltonian diffeomorphism of $T^*S^1$ lifts to a connection preserving action on $L$. This lift is defined up to an homothety  on $L$. 

Let $J\rightarrow T^*S^1$ be the frame bundle of $L$. It is by construction a 3-dimensional contact manifold, with contact form $\beta$ and equipped with an $\mathbb R$ action whose flow is a Reeb flow of the contact form. 

From the discussion of Paragraph \ref{jetspace}, we conclude that $J$ is identified 
with $J^1(S^1,\mathbb R)$ with its contact form and $\mathbb R$-action.

\subsubsection{Towards exact symplectic homeomorphisms} 
The aim of the following sections is
\begin{itemize}
\item to recall basic facts about symplectic and Hamiltonian actions on the Annulus, and in particular that every exact symplectic action lift to a connexion preserving action on $L$,
\item to extend this construction to a situation with less regularity, and in particular to give a "symplectic" interpretation of the action of 
$$
\cdiff,
$$
on $T^*S^1$.
\end{itemize}
\subsection{Exact symplectomorphisms}\label{sympham}
In this section, we recall  basic facts about symplectic  diffeomeomorphisms. Let $(M,\omega)$ be a symplectic manifold, such that $\omega$ is the curvature of an orientable real line bundle $L$ equipped with a connection $\nabla$. For any curve $\gamma$ joining $x$ and $y$ we denote by $Hol(\gamma):L_x\rightarrow L_y$ the holonomy of the connection $\nabla$ along $\gamma$. If $\gamma$ is a closed curve ({\it i.e:} $x=y$), we identify $GL(L_x)$ with $\mathbb R$ and consider $Hol(\gamma)$ as a real number.

Let  $\sigma$ be a nonzero section  of $L$. Let $\beta$ be the primitive of $\omega$ defined by
$$\nabla_X\sigma=\beta(X)\sigma.$$
We observe that if $\gamma$ is a smooth closed curve, then
$$
Hol(\gamma)=e^{\int_{\gamma}\beta}.
$$
\subsubsection{Exact symplectic diffeomorphisms}
A symplectic  diffeomorphism is of $M$ is {\em exact} if for any closed curve 
\begin{eqnarray}
Hol(\sigma)=Hol(\phi(\sigma)).\label{holcond}
\end{eqnarray}
The following Proposition clarifies this last condition
\begin{proposition}\label{holinv} Let $\phi$ be a symplectic diffeormorphism. The following conditions are equivalent \begin{itemize}
\item for any curve $\sigma$, 
$
Hol(\sigma)=Hol(\phi(\sigma))$.
\item $\phi^*\beta -\beta$ is exact,
\item The action of $\phi$ lifts to a connection preserving action on $L$.
\end{itemize}
Furthermore if $M=T^*S^1$ and $\sigma_0$ is the zero section of $T^*S^1\rightarrow S^1$, then
the map
$$
\gamma:\mapping{H}{\mathbb R}{\phi}{Hol(\phi(\sigma_0))}.
$$
is a group homomorphism, whose kernel is the group of exact diffeomorphisms
\end{proposition}
\proof
The only point which is not obvious is the fact that for an exact symplectic diffeomorphism $\phi$, the action lift to an action on the  line bundle $L$.
We define a lift $\hat\phi$ in the following way. We first choose a base point $x_0$ in $T^*S^1$ and $\eta$ a curve joining $x_0$ to $\phi(x_0)$. For any point $x$ in $T^*S^1$, we choose a  curve $\gamma$ joining $x_0$ to $x$, then we define
$$
\hat\phi_x= Hol(\phi(\gamma))Hol(\eta)Hol(\gamma)^{-1}:L_x\rightarrow L_{\phi(x)}.
$$
Since for every closed curve
$$
Hol(\sigma)=Hol(\phi(\sigma)),
$$
the linear map $\hat\phi_x$ is independent of the choice of $\gamma$.
Finally we define
$$
\hat\phi(x,u)=(\phi(x),\hat\phi_x u).
$$
By construction, $\hat\phi$ preserve the parrallel transport along any  curve $\gamma$, that is
$$
\hat\phi(Hol(\gamma)u)=Hol(\phi(\gamma)).\hat\phi(u).
$$
\qed
\subsection{$\pi$-Symplectic Homeomorphism}
The action of $Di\!f\!f(S^1)$ on $S^1$ lifts to an action by area preserving homeomorphisms on $T^*S^1=T^*S^1$. On the other hand $\Omega^1(S^1)$, the vector space of continuous one-forms on $S^1$ also acts on $T^*S^1$ by area preserving homeomorphism in the following way 
$$
fd\theta.(\theta, t)=(\theta,t+f(\theta)).
$$
We observe that $Di\!f\!f(S^1)$ normalizes this action.

We define the group of {\em $\pi$-symplectic homeomorphisms} to be the group $\Omega^1(S^1)\rtimes Di\!f\!f(S^1)$. Here is an immediate characterisation of  $\pi$-symplectic homeomorphisms whose proof is identical as that of Lemma \ref{areapreserve}

\begin{proposition}\label{pisympl}
An area preserving homeomorphism $\phi$ of $T^*S^1$ is  $\pi$-symplectic, if there is $\psi$ in $Di\!f\!f(S^1)$, such that
$$
\pi\circ\phi=\psi\circ\pi.
$$
\end{proposition}
\subsubsection{$\pi$-curves, holonomy of $\pi$-curves}

A {\em $\pi$-curve} is a continuous curve $c=(c_\theta, c_r)$ with values in $T^*S^1$ such that $c_\theta$ is $C^1$. We define
$$
\int_c \beta = \int c_r dc_\theta,
$$
and 
$$
Hol(c)=e^{\int_c\beta}.
$$
We collect in the following Proposition a few elementary facts.
\begin{proposition}
\begin{enumerate}
\item The image of a $\pi$-curve by a $\pi$-symplectic homeomorphism is a $\pi$-curve,
\item Let $\Phi=(\alpha,\phi)$ be a $\pi$-symplectic homeomorphism.
Let $c$ be a closed $\pi$-curve, then
$$
Hol(\Phi(c))=Hol(c)e^{\int_{S^1} \alpha}.
$$
\end{enumerate}
\end{proposition}
\subsubsection{$\pi$-exact symplectomorphism}
We finally define the group of {\em $\pi$-exact symplectomorphism} to be the group 
$$
Ex^\pi = (C^1(S^1)/\mathbb R)\rtimes Di\!f\!f(S^1),
$$
where $C^1(S^1)/\mathbb R$ is identified with the space of continuous exact one-forms on $S^1$.
Similarly, we finally  the group of {\em $\pi$-Hamiltonian} to be the group 
$$
Ham^\pi = (C^1(S^1)/\mathbb R)\rtimes Di\!f\!f_+(S^1),
$$
By the above remarks, we deduce immediately, reproducing the proof of Proposition \ref{holinv}
\begin{proposition}\label{holeq}
\begin{enumerate}
\item The group of $\pi$-exact symplectomorphisms is the group of $\pi$-symplectic homomorphisms $\phi$ such that for any closed $\pi$-curve $c$,
$$
Hol(c) =Hol(\phi(c)).
$$
\item The action of any $\pi$-exact symplectomorphism $\Phi$ lift to an action of a homeomorphism $\hat\Phi$ on $L$ such that for any $\pi$-curve
\begin{equation}
\hat\Phi(Hol(c))=Hol(\phi(c))\hat\Phi \label{holinvpi}.
\end{equation}
\item Furthermore, $\hat\Phi$ is determined by Equation (\ref{holinvpi}) up to a multiplicative constant. As a consequence there exists a  element $(f,\phi)$ of $C^1(S^1)\rtimes Di\!f\!f(S^1)$ such that the induced action of  $\hat\Phi$ on $J$, identified with the frame bundle of $L$ is $(f,\phi)$. We observe that $f$ is defined uniquely up to an additive constant.
\item Finally, if $\Phi$ is Hölder, so is $\hat\Phi$.
\end{enumerate}
\end{proposition}

From this we obtain the following characterisation of $\cdiff$ which generalises Proposition \ref{cidiff}.

\begin{coro}\label{caractcdiff}
A Hölder homeomorphism of $J$ belongs to $\cdiff$ if and only if 
\begin{enumerate}
\item it preserves $\mathcal F$, 
\item it commutes with the canonical flow,
\item it is above a $\pi$-exact symplectomorphism of $T^*S^1$.
\end{enumerate}
\end{coro}
\subsubsection{Width}
We generalise the notion of width  for a $\pi$-exact symplectomorphism $\phi$ in the following obvious way.
\begin{definition}{\sc [Width, Action difference]}
Let $x$ and $y$ two fixed points of $\phi$. Let $c$ be a $\pi$-curve joining $x$ to $y$. We define the {\em action difference} of $x$ and $y$ by
$$
\delta(\phi;x,y)=Hol(c\cup\phi(c))=e^{\int_c\beta-\int_\phi(c)\beta}.
$$
By Proposition \ref{holeq},it follows that this quantity does not depend on $c$. Furthermore
$$
\delta(\phi;x,y)=\hat\phi(x)/\hat\phi(y).
$$
Finally, the {\em width} of $\phi$ is
$$
w(\phi)=\sup_{x,y}\delta(\phi;x,y).
$$
\end{definition}
This quantity is invariant by conjugation under $\pi$-symplectic homeomorphisms. We  identify this with the period of some cross ratio.
 
\begin{proposition} \label{width}Let $\rho$ be a representation in $\Hom_H$.  Let $l_\rho$ be the spectrum of $\rho$ ({\it cf} Paragraph \ref{crosshomh}) and $b_\rho$ be its associated cross ratio  ({\it cf} Paragraph \ref{specthomh}) with periods $l_{b_\rho}$. Then
$$
\log(w(\rho(\gamma))= l_\rho(\gamma)+l_\rho(\gamma^{-1})=2l_{b_\rho}(\gamma).
$$
\end{proposition}

\proof Indeed, from the definition of the spectrum, for the attractive fixed point $\gamma^+$ of $\rho(\gamma)$ we have 
$$
\hat\phi(\gamma^+)=e^{l_\rho(\gamma)}.
$$
Since the repulsive point of $\rho(\gamma)$ is the attractive fixed point of $\rho(\gamma^{-1})$, we obtain that
$$
\log(w(\rho(\gamma))= l_\rho(\gamma)+l_\rho(\gamma^{-1}).
$$
We now recall that according to Proposition \ref{conjug1cross}, there exists a $\grf$-equivariant map $\Theta$ from $T^*S^1$ to $\bgrf^{2*}$ which sends the fibres of $T^*S^1\rightarrow S^1$ to the first factor. We define the inverse images of the second factor as {\em horizontal} curves. If two points $a$ and $b$ belong to the same fibre of $T^*S^1$, $[a,b]$ is the arc joining $a$ and $b$ long that fibre
 
Let now $x\in T_{\gamma^+}S^1$ and $y\in T_{\gamma^-}S^1$ be the two fixed points of $\rho(\gamma)$ in $T^*S^1$.  Let $t$ be a point in $T_{\gamma^+}S^1$. Let $u$ be a point in $T_{\gamma^+}S^1$ which is joined to $t$ by an horizontal arc $c$. 

We note that by definition
$$
l_{b_\rho}(\gamma)=\frac{1}{2}\log(Hol(c\cup [t,\rho(\gamma)(t)]\cup \gamma(c)\cup [\rho(\gamma)(u),u])).
$$
Let now $\tilde c = [x,t]\cup c\cup [u,y]$. We then have
$$
l_{b_\rho}(\gamma)=\log(Hol(\tilde c \cup \rho(\gamma)(\tilde c))=\log(w(\rho(\gamma))).
$$
The Assertion follows
 \qed

\subsection{Conjugation Lemma}
We state the hypothesis and notations of the result we prove in this paragraph
\begin{itemize}
\item Let $\kappa$ be a a Hölder homeomorphism of $S^1$. Let 
$$
D_{\kappa}=S^1\times S^1\setminus\{(s,t)/\kappa(s)=t\}.
$$
\end{itemize}
We define
\begin{itemize}
\item Let $\pi_D$ be the projection from $D_{\kappa}$ on the first $S^1$ factor.
\item Let $ds$ (resp.  $dt$) be Lebesgue probability measure on the first (resp. second) factor of $S^1\times S^1$.
\item Finally, let $d\theta\otimes dr$ be the measure associated to the canonical symplectic form of $T^*S^1$.
\end{itemize}
We also consider
\begin{itemize}
\item Let $f$ be a positive continuous Hölder function such that
\begin{eqnarray}
\forall s,\ t \hbox{ such that }\kappa(s)\not=t,  \int_t^{\kappa(s)}f(s,u)du=  \int_{\kappa(s)}^t f(s,u)du=\infty.\label{infmesure}
\end{eqnarray}
\end{itemize}
Finally,
\begin{itemize}
\item Let $\omega=f(s,t)ds\otimes dt$. Let
$$
\beta(s,t)=\bigg(\int_{g(s)}^t f(s,u)du\bigg) ds.
$$
Note that $\beta$ is  a primitive of $\omega$. 
\item For any $C^1$-closed curve $c$, let $$
Hol_\omega(c)=e^{\int_c\beta}.
$$
\end{itemize}
\begin{lemma}\label{conjug}Assume that $f$ satisfies  Hypothesis (\ref{infmesure}). Then there exists a Hölder homeomorphism $\psi$ from $D_\kappa$ to $T^*S^1$ such that
\begin{itemize}
\item $\psi_* \omega=d\theta\otimes dr$,
\item $\pi \circ \psi=\pi_D$. 
\item  if $\gamma_1$, $\gamma_2$ are two $C^1$ diffeomorphisms with Hölder derivatives of  the circle such that $\gamma=(\gamma_1,\gamma_2)$ preserves $\omega$ then
$$
\psi\circ\gamma\circ\psi^{-1},
$$
is Hölder and is $\pi$-symplectic above $\gamma_1$,
\item if $c$ is a $C^1$-curve in $D$, then $\psi(c)$ is a $\pi$-curve. Furthermore if $c$ is a closed curve then $Hol_\omega(c)=Hol(\psi(c))$.
\end{itemize}
The homeomorphism $\psi$ is unique up to right composition with a $\pi$-exact symplectomorphism. Finally $\psi$ depends continuously  on $\kappa$ and $f$. \end{lemma}
\proof
The uniqueness part of this  statement follows by the characterisation of $\pi$-exact hamiltonian given in  Proposition \ref{pisympl}.
The proof is completely explicit. Let $g$ be a $C^1$-diffeomorphism such that
$\forall s, g(s)\not=\kappa(s)$ 
We consider 
$$
\psi :
\mapping{D}{T^*S^1=S^1\times \mathbb R}{(s,t)}{(s,\int_{g(s)}^t f(s,u)du)}.
$$
It is immediate to check that
\begin{itemize}
\item $\psi_* \omega=d\theta\otimes dr$,
\item $\pi\circ \psi=\pi_D$. 
\end{itemize}
We also observe that $\psi$ is a  Hölder map and a homeomorphism. We prove now that $\psi^{-1}$ is Hölder. We have
$$
\psi^{-1}(s,u)=(s,\alpha(s,u)).
$$ 
Where
$$
\int_{g(s)}^{\alpha(s,u)}f(s,w)dw=u.
$$
It is enough to prove that $\alpha$ is Hölder. We work locally, so that $f$ is bounded from below by $k$. Firstly,
$$
\vert u- v\vert=\vert \int_{\alpha(s,v)}^{\alpha(s,u)}f(s,w)dw\vert\geq k\vert\alpha(s,v)-\alpha(s,u)\vert.$$
This prove that $\alpha$ is Hölder with respect to the second variable. Let us take care of the first variable. We have 
$$
\int_{g(s)}^{\alpha(s,u)}f(s,w)dw=u=
\int_{g(t)}^{\alpha(t,u)}f(t,w)dw.
$$
Hence,
$$
\int_{\alpha(s,u)}^{\alpha(t,u)}f(t,w)dw= \int_{g(t)}^{g(s)}f(t,w)dw +\int_{g(t)}^{\alpha(t,u)}(f(t,w)-f(s,w))dw.
$$
Since we work locally, we may assume that $k\leq\vert f(t,w)\vert\leq K $ and 
$$
\vert f(t,w)-f(s,w)\vert\leq C\vert t-s\vert^\beta,
$$
we have
$$
k\vert \alpha(s,u)-\alpha(t,u)\vert  \leq K \vert t-s\vert + C\vert t-s\vert^\beta.
$$
Therefore $\alpha$ is a Hölder function.

By the construction of $\psi$, if $c$ is a $C^1$-curve in $D$, then $\psi(c)$ is a $\pi$-curve.
Recall that 
$$
\beta(s,t)=\bigg(\int_{g(s)}^t f(s,u)du\bigg) ds
$$
is a primitive of $\omega$. It follows that  if $c:s\mapsto (c_{1}(s),c_{2}(s))$ is a $C^1$ curve in $D_{\kappa}$. Then
$$
\log(Hol_\omega (c))=\int_c\beta= \int_{S^1}\!\!\int_{g(s)}^{c_{2}(s)} f(s,u)\dot c_{1}(s) duds=\log(Hol(\psi(c))).
$$

It remains to check that $\psi\circ\gamma\circ\psi^{-1}$ is $\pi$-symplectic and Hölder. But this is immediate by construction. The continuity of $\psi$ on $\kappa$ and $f$ follows from the construction.
\qed

\subsection{Proof of  Theorem \ref{conjugtheo}}
\subsubsection{A preliminary lemma}
 
We need the following lemma
\begin{lemma}\label{lemmconjug}
Let   $\rho_1$ and $\rho_2$ be two  $H$-Fuchsian representations of $\grf$ in $Di\!f\!f^{h}(S^1$). Let 
$\kappa$ be  a Hölder homeomorphism of $S^1$,
such that $\kappa\circ\rho_1=\rho_2\circ\kappa$. Let $f(s,t)$ be a positive continuous function on 
$$
D_{\kappa}=S^1\times S^1\setminus \{(s,t)/\kappa(s)=t\}.
$$
Assume that, for all $\gamma$ in $\grf$,  $\omega=f(s,t)ds\otimes dt$ is invariant under the action of $(\rho_1(\gamma),\rho_2(\gamma))$. Then
$$
\forall (s,t)\in D_{\kappa}, \int_{t}^{\kappa(s)}f(s,u)du=\int_{\kappa(s)}^{t}f(s,u)du=\infty.
$$
\end{lemma}
\proof
For the sake of simplicity, we write $\rho_i(\gamma)=\gamma^i$. We first observe that the invariance of $\omega$ yields
$$
f(s,t)=\frac{d\gamma^1_n}{ds}(s)\frac{d\gamma^2_n}{dt}(t)f(\gamma^1_n(s),\gamma^2_n(t)).
$$
For any $(s,t)$ in $D_\kappa$, we may find a sequence $\{\gamma_n\}_{n\in\mathbb N}$ such that
\begin{eqnarray}
\lim_{n\rightarrow\infty}\gamma^1_n (s)&=& s_0,\label{t1}\\
\lim_{n\rightarrow\infty}\gamma^2_n (t)&=&t_0\not=\kappa(s_0)\label{t2},\\
\lim_{n\rightarrow\infty}\frac{d\gamma^1_n}{ds}(s)&=&+\infty.\label{t3}
\end{eqnarray}
Now using the invariance of $\omega$ by $(\gamma^1,\gamma^2)$, we obtain that
\begin{eqnarray*}
\int_{\kappa(s)}^t f(s,u)du &=& \int_{\kappa(s)}^t  \frac{d\gamma^1_n}{ds}(s)\frac{d\gamma^2_n}{du}(u)f(\gamma^1_n(s),\gamma^2_n(u))du\cr &=&
\frac{d\gamma^1_n}{ds}(s) \int_{\kappa(s))}^{ t}\frac{d\gamma^2_n}{du}(u)f(\gamma^1_n(s),\gamma^2_n(u))du\cr
&=&
\frac{d\gamma^1_n}{ds}(s) \int_{\gamma^2_n(\kappa(s))}^{\gamma^2_n (t)}f(\gamma^1_n(s),u)du\cr
&=&
\frac{d\gamma^1_n}{ds}(s) \int_{\kappa(\gamma^1_n( s))}^{\gamma^2_n (t)}f(\gamma^1_n(s),u)du
\end{eqnarray*}
From  Assertions (\ref{t1}) and (\ref{t2}), we deduce that 
$$
\lim_{n\rightarrow\infty}\int_{\kappa(\gamma^1_n( s))}^{\gamma^2_n (t)}f(\gamma^1_n(s),u)du
 = \int_{\kappa(s_0)}^{t_0}f(s_0,u)du>0.
$$
Hence Assertion (\ref{t3}) shows that 
$$
\int_{\kappa(s)}^t f(s,u)du =\infty.
$$
A similar argument shows that 
$$
\int_{t}^{\kappa(s)} f(s,u)du =\infty.
$$
\qed
\subsubsection{Proof of Theorem \ref{conjugtheo}} 
\proof  Recall our hypothesis and notations.
Let 
$\kappa$ be  a Hölder homeomorphism of $S^1$. Let 
$$
D_{\kappa}=S^1\times S^1\setminus\{(s,t)/\kappa(s)=t\}.
$$
such that $\kappa\circ\rho_1=\rho_2\circ\kappa$. Let $p: M=\mathbb R\times D_{\kappa}\rightarrow D_{\kappa}$ be a principal $\mathbb R$-bundle over $D_{\kappa}$ equipped
with a connection $\nabla$. Assume  the  curvature $\omega$ of $\nabla$ is such that $\omega=f(s,t)ds\wedge dt$ with $f$ positive and Hölder.  . 
Let 
\begin{eqnarray*}
\pi&:&\mapping{M}{S^1}{(u,s,t)}{t}\cr
\pi_1&:&\mapping{M}{\mathbb R\times S^1}{(u,s,t)}{(u,s)}\cr
p &:& \mapping{M}{D_{\kappa}}{(u,s,t)}{(s,t)}
\end{eqnarray*} 
Assume that $\grf$ acts on $M$ by $C^1$-diffeomorphisms with Hölder derivatives. Assume this action preserves $\nabla$, and  that 
\begin{itemize}
\item $M/\grf$ is compact
\item The action of $\mathbb R$ on $M/\grf$ contracts the fibres of $\pi_{1}$.
\item There exists two $H$-Fuchsian  representations $\rho_{1}$ and $\rho_{2}$ of $\grf$, such that
$$
p(\gamma u)= (\rho_{1}(\gamma)(p(u)),\rho_{2}(\gamma)(p(u))).
$$  
\end{itemize}

We want to prove  there exists a $\mathbb R$-commuting Hölder homeomorphism $\hat\psi$ from $M$ to $J$ over a homeomorphism $\psi$ from $D$ to $T^*S^1$, a representation $\rho$ of $\grf$ in $\cdiff$ element of $\Hom_{H}$, such that

$$ \hat\psi\circ\gamma=\rho(\gamma)\circ\hat\psi.$$
Notice first that by Lemma \ref{lemmconjug}
$$
\forall (s,t)\in D_{\kappa}, \int_{t}^{\kappa(s)}f(s,u)du=\int_{\kappa(s)}^{t}f(s,u)du=\infty.
$$
It follows the hypotheses of  Lemma \ref{conjug} are satisfied. Thus, there exists a Hölder homeomorphism $\psi$ from $D$ to $T^*S^1$ unique up to right composition with a $\pi$-exact hamiltonian, such that
\begin{enumerate}
\item $\psi_* \omega=d\theta\otimes dr$,
\item $\pi \circ \psi=\pi_D$. 
\item  if $\gamma_1$, $\gamma_2$ are two $C^1$ diffeomorphisms with Hölder derivatives of  the circle such that $\gamma=(\gamma_1,\gamma_2)$ preserves $\omega$ then
$$
\psi\circ\gamma\circ\psi^{-1},
$$
is Hölder and $\pi$-symplectic above $\gamma_1$,
\item if $c$ is a $C^1$-curve in $D$, then $\psi(c)$ is a $\pi$-curve and $Hol_\omega (c)=Hol(\psi(c))$.
\end{enumerate}

Let $\gamma$ be an element of $\grf$. Since $\gamma$ acts on $M$ preserving the connection $\nabla$, it follows that $f(\gamma)=(\rho_{1}(\gamma),\rho_{2}(\gamma))$ acts on $D_{\kappa}$ in such  a way that for all $C^1$-curves
$Hol_\nabla(c)=Hol_\nabla(f(\gamma)c)$. 

We now show that 
\begin{eqnarray}
Hol_\omega(c)=Hol_\omega(f(\gamma)c).\label{connect}
\end{eqnarray} 
We choose a trivialisation of $L$. In this trivialisation $\nabla= D+\beta+\alpha$, where $D$ is the trivial connection and $\alpha$ is a  closed form. Then
$$
Hol_\nabla(c)=Hol_\omega(c).e^{\int_c\alpha}.
$$
To prove Equality (\ref{connect}), it suffices to show 
\begin{eqnarray}
\int_{f(\gamma)(c)}\alpha=\int_c\alpha.\label{connectfin}
\end{eqnarray}
But $\rho_1(\gamma)$ preserves the orientation of $S^1$, hence is connected to the identity by a family of mapping $f_t$. It follows that  $f(\gamma)$ is also homotopic to the identity through the family $(f_t, \kappa f_t\kappa^-1)$. This implies that $f(\gamma)$ acts trivially on the homology. Hence Equality (\ref{connectfin}).

It follows from (3) and (4) that $g(\gamma)=\psi\circ f(\gamma)\circ\psi^{-1}$ is a $\pi$-exact symplectomorphism.  

Finally, we define a map $\hat\psi$ from $M$ to $J$ above $\psi$, commuting with $\mathbb R$ action, ''preserving the holonomy '' well defined up right composition by an element of $\mathbb R$. Let fix an element $y_{0}$ of the fibre of $p$ above $x_{0}$ in $D_{\kappa}$ and an element $z_{0}$ of the fibre of $\pi_{2}$ above $\psi(x_{0})$.  Let $y$ be an element of the fibre of $p$ above some point $x$ of $D_{\kappa}$. Let $c$ be path joining  $x_{0}$ to $x$. We observe that there exists $\mu$ in $\mathbb R$ such that
$$
y=\mu + Hol_\nabla(c)y_0.
$$
We define 
$$
\hat\psi(y)=\mu + Hol(\psi(c))z_{0}.
$$
Since for any closed curve $Hol_\nabla (c)=Hol(\psi(c))$, it follows that $\hat\psi(y)$ is independent of the choice of $c$.

By construction, $\hat\psi$ satisfies the required conditions. The uniqueness statement follows by the Corollary \ref{caractcdiff} The continuity statement follows by the corresponding continuity statement of Lemma \ref{conjug} and the construction. \qed 

\section{Negatively curved metrics}\label{negcurvsect}
 We first use our conjugation Theorem \ref{conjugtheo} to prove the space $\Rep_{H}$ contains an interesting space.
\begin{theorem}\label{negcurv}
Let $\mathcal M$ be the space of negatively curved metrics on the surface $S$. Then there exists a continuous injective map $\psi$ from $\mathcal M$ to $\Rep_{H}$. Furthermore, for any $\gamma$ in $\grf$
$$
l_{g}(\gamma)=l_{\psi(g)}(\gamma).
$$
Here $l_{g}(\gamma)$ is the length of the closed geodesic for $g$ freely homotopic to $\gamma$, and $l_{\psi(g)}$ is the $\psi(g)$
-length of $\gamma$.
Finally, $\psi(g_{0)}=\psi(g_{1})$, if and only if there exists a diffeomorphism $F$of $S$,  homotopic to the identity, such that $F^*(g_{0})=g_{1}$. 
\end{theorem}

We first recall some facts about the geodesic flow and the boundary at infinity of negatively curved manifolds.

\subsection{The boundary at infinity and the geodesic flow}

Let $S$ be a compact surface equipped with a negatively curved metric. Let $\tilde S$ be its universal cover. Let $U\tilde S$ (resp. $US$) be the unitary tangent bundle of $\tilde S$ (resp. $S$). Let $\phi_{t}$ be the geodesic flow on these bundles. Let $\partial_{\infty}\tilde S$ be the boundary at infinity of $\tilde S$.

We collect in the following Proposition well known facts :  

\begin{proposition}\label{factnegcurv1}
\begin{enumerate}
\item $\partial_{\infty}\tilde S=\bgrf$.
\item $\partial_{\infty}\tilde S$ has a $C^1$-structure depending on the choice of the metric such that the action of $\grf$ on it is by $C^{1,h}$-diffeomorphisms. The action of $\grf$ is Hölder conjugate to a Fuchsian one.
\item Let $\mathcal G$ be the space of geodesics of the universal cover of $S$. Then $\mathcal G$ is $C^{1,h}$-diffeomorphic to $(\bgrf)^{2*}$.
\item $U\tilde S \rightarrow \mathcal G$, is a $\mathbb R$-principal bundle (with the action of the geodesic flow). Furthermore the Liouville form is a connection form for this action invariant under $\grf$, and its  curvature is symplectic.
\end{enumerate}
\end{proposition}

The following well known Proposition explains the identification of the unitary tangent bundle with a suitable product of the boundary at infinity.

\begin{proposition}\label{factnegcurv2}
There exist a homeomorphism $f$ of 
$$
\bgrf^{3*}=\{ \hbox{ oriented } (x,y,z)\in(\bgrf)^{3}/ x\not=y\not= z\not= x\}
$$
with $U\tilde S$
such that if $\psi_{t}=f^{-1}\circ\phi_{t}\circ f$ then, 
\begin{itemize}
\item $\psi_{t}(z,x,y)=(w,x,y)$,
\item let $\seq{t}$ be a sequence of real numbers going to infinity, let $(x,y,z)\in\bgrf^{3*}$, let $\seq{\gamma}$ be a sequence of elements of $\grf$, such that 
$$
\{\gamma_{m}\circ\psi_{t_{m}}(z,x,y)\}\inn \hbox{ converges to } (z_0,x_{0},y_{0}).
$$ 
Then for any $w$, $v$ such that $(w,x,v)\in\bgrf^{3*}$, there exists $v_{0}$ such that
$$
\{\gamma_{m}\circ\psi_{t_{m}}(w,x,v)\}\inn\hbox{ converges to } (v_{0},x_{0},y_{0}).
$$ 
\end{itemize}
\end{proposition}

\proof We explain first the construction of the map $f$. Let $(z,x,y)\in\bgrf^{3*}$. Let $\gamma$ the geodesic in $\tilde S$ going from $x$ to $y$. Let $\gamma(t_{0})$ be the projection of $y$ on $\gamma$, that is the unique minimum on $\gamma$ of the horospherical function associated to $z$. We set $f(z,x,y)=\dot\gamma(t_{0})$. Then the first property of $f$ is obvious, and the second one is a classical consequence of negative curvature, namely that two geodesics with the same endpoints at infinity go exponentially closer and closer. \qed
 
\subsection{Proof of Theorem \ref{negcurv}}

We first have to check the hypotheses of Theorem \ref{conjugtheo} are satisfied. This follows at once by Proposition \ref{factnegcurv1}  and \ref{factnegcurv2}, using Proposition \ref{conditionalternate}. Therefore associated to a negatively curved metric we obtain a representation in $\Hom_*$. For a hyperbolic mtric, we obtain precisely a $\infty$-Fuchsian representation. Since the space of negatively curved metric is connected, all the representation are actually in $\Hom_H$.

Finally if $\psi(g_0)=\psi(g_1)$, then the two metrics have the same length spectrum and  therefore are isometric by Otal's Theorem \cite{JPO}. This proves injectivity. The continuity follows by the continuity statement of Theorem \ref{conjugtheo}.

\section{Hitchin component}\label{hichsect}

Let $\Rep_{H}(\grf,SL(n,\mathbb R))$ be a {\em Hitchin component}, {\em i.e.} a connected component of the space of representations of $\grf$ in 
$SL(n,\mathbb R)$ containing the Fuchsian representations as defined in \cite{H} ( see also \cite{FL3}).

We now prove that $\Rep_H$ contains all these Hitchin components.

\begin{theorem}\label{hitchcross}
There exists a continuous injective  map 
$$
\psi : \Rep_{H}(\grf,SL(n,\mathbb R))\rightarrow \Rep_{H},
$$
such that, if $\rho\in \Rep_{H}(\grf,SL(n,\mathbb R))$ then
\begin{itemize}
\item for any $\gamma$ in $\grf$, we have 
\begin{equation}
l_{\psi(\rho)}(\gamma)=\log(\vert\frac{\lambda_{max}(\rho(\gamma)}{\lambda_{min}(\rho(\gamma)}\vert).\label{spect=vp}
\end{equation}
Here, $l_{\psi(\rho)}(\gamma)$ is the $\psi(\rho)$-length of $\gamma$, and $\lambda_{max} (a)$ (resp. $\lambda_{min}(a)$) denote the maximum (resp. minimum) real eigenvalue of the endomorphism $a$ in absolute value.
\item The cross ratio associated to $\rho$ and $\psi(\rho)$ coincide.
\item $\psi(\rho)$ is compatible ({\it cf} Definition \ref{compa}).
\end{itemize}
\end{theorem}

In some sense, this result says that $\cdiff$ is a version of $SL(\infty,\mathbb R)$.

\subsection{Hyperconvex curves}\label{hypercurvehitch}
For every Hitchin representation $\rho$ in $\sln$, there exists a $\rho$-equivariant hyperconvex map $\xi$ from $\bgrf$ to $\pnu$. More precisely, in \cite {FL3} we proved there exists a $\rho$-equivariant Hölder map $(\xi^{1},\xi^{2},\ldots,\xi^{n-1})$ from $\bgrf$ to the flag manifold $\mathcal F$, called the {\em osculating flag of $\xi$}  or the {\em limit curve of $\rho$},
such that 
\begin{itemize} 
\item $\xi=\xi^{1}$,
\item $\xi^{p}$
  is with values in the Grassmannian of $p$-planes, 
  \item if
  $(n_1,\ldots,n_l)$ are positive integers such that
  $$
  \sum_{i=1}^{i=l}n_i\leq n, 
  $$ 
  if $(x_{1},\ldots,x_{l})$ are
  distinct points, then the following sum is direct
\begin{eqnarray}
\xi^{n_i}(x_i)\oplus\ldots\oplus\xi^{n_{l}}(x_{l});\label{fre2}
\end{eqnarray}
\item finally, for every $x$, let $p=n_{1}+\ldots+n_{l}$,
  then
\begin{eqnarray}
\lim_{(y_1,\ldots,y_l)\rightarrow x, y_i
{\hbox{\tiny all distinct}}}
(\bigoplus_{i=1}^{i=l}\xi^{n_i}(y_i))=\xi^{p}(x).\label{fre3}
\end{eqnarray}
\end{itemize}
It follows in particular that $\xi^1(\bgrf)$ (resp. $\xi^{n-1}(\bgrf)$) is a $C^1$-curve with Hölder derivatives in $\pnu$ (resp.$\pnud$). Furthermore, for $x\not=y$, the sum $\xi^1(x)\oplus \xi^{n-1}(y)$ is direct.

We may rephrase these properties in the following way

\begin{proposition}\label{hitch2}Let $\rho$ be  a hyperconvex representation. Let $\xi=(\xi^1,\ldots,\xi^{n-1})$ be the limit curve of $\rho$. Then, 
there exist
\begin{itemize}
\item Two $C^1$ embeddings with Hölder derivatives, $\eta_1$ and $\eta_2$, of $S^1$ in respectively $\pnu$ and $\pnud$,
\item Two representations $\rho_1$ and $\rho_2$ of $\grf$ in $Di\!f\!f^{h}(S^1$), the group of $C^1$-diffeomorphisms of $S^1$ with Hölder derivatives,
\item a Hölder homeomorphism $\kappa$ of $S^1$,
\end{itemize}
such that
\begin{enumerate}
\item $\eta_1(S^1)=\xi^1(\bgrf)$ and $\eta_2(S^1)=\xi^{n-1}(\bgrf)$,
\item $\eta_i\circ\rho_i=\rho\circ\eta_i$,
\item if $\kappa(s)\not=t$, then the sum $\eta_1(s)\oplus\eta_2(t)$ is direct.
\item $\kappa\circ\rho_1=\rho_2\circ\kappa$,
\end{enumerate}
\end{proposition}

\proof Let $\eta_1$ (resp. $\eta_2$) be the arc-length parametrisationof $\xi^1(\bgrf)$ (resp. $\xi^{n-1}(\bgrf)$). Let  
$$
\kappa= (\eta_2)^{-1}\circ\xi^1\circ(\xi^{n-1})^{-1}\circ\eta_1.$$
The result follows.
\qed

For later use we introduce the continuous map
\begin{equation}
\dot\eta=(\eta_1^{-1}\circ\xi^1,\eta_2^{-1}\circ\xi^{n-1}),\label{doteta}
\end{equation}
from $\bgrf^{2*}$ to $D_\kappa$.

\subsection{Proof of Theorem \ref{hitchcross}}
We use in this section the independent results proved in the Appendix \ref{sympnat}.

Let $\rho$ be a hyperconvex representation. Let $\eta_1$, $\eta_2$ and $\kappa$ as in Proposition \ref{hitch2}. Let $\eta=(\eta_1,\eta_2)$. Let  
$$
D_\kappa=(S^1)^2\setminus\{(s,t)/\kappa(s)=t\},
$$
We observe that $\grf$ acts by $\dot\rho=(\rho_1,\rho_2)$ on $D_\kappa$.
It follows by Proposition \ref{hitch2}, that $\eta$ is a $\rho$-equivariant (with respect to the action on $D_\kappa$ given by $\dot\rho$) $C^1$ map
to
$$
\mathbb P(n)^{2*}=\pnu\times\pnud\setminus\{(D,P)/D\subset P\}.
$$
In  Section \ref{sympproj} of the Appendix, we show  there exists a $SL(n,\mathbb R)$ invariant $\mathbb R$-principal  bundle $L$ on $\mathbb P (n)^{2*}$ equipped with a connection whose curvature is a symplectic form $\Omega$. 

\begin{itemize}
\item Let $M$ be the induced bundle by $\eta$: $M=\eta^*L$. It is a $\mathbb R$-principal bundle over $D_\kappa$ equipped with an action of $\grf$ coming from the $SL(n,\mathbb R)$ action on $L$. 
\item Let $\phi_t$ be the flow of the induced $\mathbb R$-action on $M$.
\item Let $\mathcal L$ be the foliation of $M$ by parallel vector fields along the curves $c_s: t\mapsto(s,t)$ in $D_\kappa$.  
\end{itemize}
We observe the following

\begin{proposition}{\label{action}}
Let $\gamma$ be an element of $\grf$. Let $x={\dot\eta}(\gamma^+,\gamma^-)$ be a fixed point of $\gamma$ in $D_\kappa$. Let $\lambda_{max}$ and $\lambda_{min}$ be the largest and smallest eigenvalues (in absolute values)  of $\rho(\gamma)$. Then the action of $\gamma$ on the fibre $M_x$ of $M$ above $x$ is given by the translation by $$
\log\big\vert\frac{\lambda_{max}}{\lambda_{min}}\big\vert.
$$
\end{proposition}
\proof This follows from the  last point of Proposition \ref{sympproj}.\qed

To complete the proof of Theorem \ref{hitchcross}, we prove
\begin{itemize}
\item $\omega= \eta^*\Omega$ is symplectic and Hölder (in Proposition \ref{symphold}),
\item $M/\grf$ is compact (in Proposition \ref{grfcomp})
\item The action of $\mathbb R$ on $M$ contracts the leaves of $\mathcal L$ (in Proposition \ref{contractleaves}).
\end{itemize}

We explain now how these properties imply our Theorem :
by Theorem \ref{conjugtheo},  there exists a homeomorphism $\hat\psi$ from $M$ to $J$, an area preserving  homeomorphism $\psi$ from $D$ to $T^*S^1$ such that
\begin{itemize}
\item $\hat\psi$ commutes with the $\mathbb R$-action,
\item
$
\pi_2\circ\hat\psi=\psi\circ p,
$ 
\item $\hat\psi$ sends $\mathcal L$ to $\mathcal F$,
\end{itemize}
and a  representation $\rho$ of $\grf$ in $\cdiff$ such that
$$ \hat\psi\circ\gamma=\rho(\gamma)\circ\hat\psi.$$
In particular, $\rho$ belongs to $\Hom^*$.

Let us prove that $\rho$ is actually in $\Hom_H$. From the definition of Hitchin component, and the fact that we can choose $\hat\psi$ to depend continuously on our parameters, it suffices to check that  for any  $n$-Fuchsian representation, {\it i.e.} a representation that factors through a Fuchsian representation in $PSL(2,\mathbb R)$ and the irreducible representation $\iota$ of this latter group in $PSL(n,\mathbb R)$. But in this case, the action of $\grf$ extends to a transitive action of $PSL(2,\mathbb R)$ and is by definition in $\Hom_H$.

The statement (\ref{spect=vp}) about the spectrum follows by Proposition \ref{action}. In particular the spectrum is symmetric. Therefore, by definition, $\psi(\rho)$ is symmetric. Finally, by construction and Proposition \ref{sympproj},  the two cross ratios coincide.

By Theorem \ref{act5},  a $\infty$-Hitchin representation which is symmetric is determined by its cross ratio. It follows that $\psi$ is injective.

Finally the continuity statement follows by the continuity statement of Theorem \ref{conjugtheo}

\subsubsection{Symplectic form}
We first prove the following
\begin{proposition}\label{symphold}
The two-form $\omega= \eta^*\Omega$ is symplectic and Hölder.
\end{proposition}
\proof The regularity of $\eta$ implies $\eta^*\Omega$ is Hölder. We begin by an observation. Let $D$ be a line in $\mathbb R^n$, $P$ a line in $\mathbb R^{n*}$, such that $D\oplus P=\mathbb R^n$.
Let $W$ be a two-plane containing $D$. Let 
$$
\hat W=T_D\mathbb P (W)\subset T_D \pnu.
$$
Let $V$ a $n-2$-plane contained in $P$. Let
$$
\hat V=T_P\mathbb P (W^{\perp})\subset T_P \pnud.
$$
From the definition of $\Omega$, we check that if 
$$
V\oplus W=\mathbb R^n,
$$
Then 
$$
\Omega\vert_{\hat V\oplus \hat W}\not=0.
$$
In the case of hyperconvex curves, 
$$
T_{\xi^1(x)}\xi^1( S^\infty)=\widehat{\xi^2(x)},\ \ \ \ \ T_{\xi^{n-1}(x)}\xi^{n-1}( S^\infty)=\widehat{\xi^{n-2}(x)}.
$$
Since by hyperconvexity $\xi^2(x)\oplus \xi^{n-2}(y)=\mathbb R^n$ for $x\not=y$, we conclude that $\eta^*\omega$ is symplectic. \qed

As a corollary, we get

\begin{proposition}\label{b=1}
Let $\rho$ be a hyperconvex representation, then $b_\rho$ is strict.
\end{proposition}
\proof Let $(x,y,z,t)$ be  a quadruple of pairwise distinct points. Let $\dot\eta$ defined as in (\ref{doteta}).
 Let $Q$ be the square in $\bgrf^{2*}$ whose vertices are $(x,y),(z,y),(x,t),(z,t)$. We know that
$$
\vert b_\rho(x,y,z,t)\vert =e^{\frac{1}{2}\int_{\dot\eta (Q)} \omega}.
$$
Since $\omega$ is symplectic, and $Q$ has a nonempty interior, we have $\int_{\dot\eta (Q)}\omega \not=0$. Hence $b_\rho(x,y,z,t)\not=1$. \qed

\subsubsection{Compact quotient}

We consider now the $\mathbb R$-principal bundle $M$ over $D_\kappa$, defined by $M=\eta^*L$. $M$ comes equipped with a connection, whose curvature form is $\omega$. Furthermore $\grf$ acts on $M$, by the pull back of the action of $SL(n,\mathbb R)$ on $L$. We now prove
\begin{proposition}\label{grfcomp}
The quotient $M/\grf$ is compact.
\end{proposition}
Recall that $\grf$ acts with a compact quotient on 
$$
\bgrf^{3*}=\{{\hbox{oriented triples }} (x,y,z)\in (\bgrf)^3 ,x\not=y, y\not= z,x\not=z\},
$$ 
We first prove :

\begin{proposition}\label{mgcompact}
There exists a continuous onto $\grf$-equivariant proper map $l$  from $\bgrf^{3*}$ to $M$. Moreover, the map $l$ is above the  map $\dot\eta=(\dot\eta_1,\dot\eta_2)$ from $\bgrf^{2*}$ to $D_\kappa$, {\em i.e} $l(z,x,y)\in M_{\dot\eta(x,y)}$.
\end{proposition}

\proof Let $(z,x,y)$ be an element of $S^{3*}=\bgrf^{3*}$.
The two transverse flags $\xi(x)$ and $\xi(y)$ defines a decomposition
$$
\mathbb R^n =L_1(x,y)\oplus L_2(x,y)\oplus\ldots\oplus L_n(x,y),
$$
such that $\xi^1(x)=L_1(x,y)$ and 
\begin{equation}
\xi^{n-1}(y)=L_{2}(x,y)\oplus \ldots\oplus L_{n}(x,y).\label{lxy}
\end{equation} 
Let $u$ be a nonzero element of $\xi^1(z)$. Let $u_i$ be the projection of $u$ on $L_i(x,y)$. By hyperconvexity, $u_i\not=0$. We choose $u$, up to sign, so that 
$$
\vert u_1\wedge\ldots \wedge u_n\vert =1.
$$
Finally we choose $f$ in $\xi^{n-1}(y)^\perp$ so that $\langle f,u_1\rangle=1$. The pair $(u_1,f)$ is well defined up to sign, and hence defines a unique element $l(z,x,y)$ of the fibre of $M$ above $\dot\eta(x,y)$. 

We prove that 
 $l$ :
$$
\mapping{S^{3*}}{M}{(z,x,y)}{l(z,x,y),}
$$
is proper. Let $\seq{(z_m,x_m,y_m)}$ be a sequence of elements of $S^{3*}$ such that $\seq{l(z_m,x_m,y_m)=(u_m,f_m)}$ converges to $(u_0,f_0)$ with $\langle f_0,u_0\rangle=1$. In particular, we have that $\seq{(x_m,y_m)}$ converges to
$(x_0,y_0)$, with $x_0\not=y_0$. We may assume after extracting a subsequence that $\seq{z_m}$ converges to $z_0$. To prove $l$ is proper, it suffices to show that
$x_0\not=z_0$ and $y_0\not=z_0$. 

Assume not, and let us suppose first that $z_0=x_0$. Let $\pi_m$ be the projection of $\xi^1(z_m)$ on $\xi^1(x_m)$ along  $\xi^{n-1}(y_m)$. We observe that $\pi_m$ converges to the identity from $\xi^1(z_0)=\xi^1(x_0)$ to  $\xi^1(x_0)$. Let $v_m\in\xi^1(z_n)$ such that $\pi_m(v_m)=u_m$. Since $\seq{u_m}$ converges to a nonzero element $u_0$, $\seq{v_m}$ converges to $u_0$. 
As a consequence, all the projections $\seq{v^i_m}$ of $v_m$ on $L_i(x_m,y_m)$ converge to zero for $i>1$. Hence,
$$
1=\vert v^1_m\wedge\ldots\wedge v^n_m\vert \rightarrow 0,
$$
and the contradiction.

Suppose now that $z_0=y_0$. Using the volume form of $\mathbb R^n$, we identify  $\xi^{n-1}(y)^\perp$ with
$$
L_{2}(x,y)\wedge\ldots\wedge L_{n}(x,y).
$$
We use the same notations as in the previous paragraph. Then $v^2_m\wedge \ldots\wedge v^n_m$ is identified with $f_m$. It follows that 
\begin{equation}
v^2_m\wedge\ldots v^n_m\rightarrow f_0.\label{limv2}
\end{equation}
Recall that by hyperconvexity
$$
\xi^1(z_m)\oplus\xi^p(y_m)\rightarrow \xi^{p+1}(y_0).
$$
Recall also that,
$$
\xi^{p}(y)=L_{n-p+1}(x,y)\oplus \ldots\oplus L_{n}(x,y).
$$
Since $\xi^1(z_m)\rightarrow \xi^1(y_0)$, it follows that for all $k$
$$
\frac{\Vert v^{1}_m\Vert}{\Vert v^{k}_m\Vert}\rightarrow 0.
$$
In particular
$$
\frac{\Vert v^{1}_m\Vert^{n-1}}{\Vert v^{2}_m\wedge\ldots\wedge v^{n}_m\Vert}\rightarrow 0.
$$
Thanks to Assertion (\ref{limv2}), we finally get that $\Vert v^{1}_m\Vert\rightarrow 0$.
Hence 
\begin{eqnarray*}
1&=&\Vert v^1_m\wedge\ldots\wedge v^n_m\Vert\\
&\leq& \Vert v^1_m\Vert \Vert v^2_m\wedge\ldots v^n_m \Vert\rightarrow 0,
\end{eqnarray*}
and the contradiction. Therefore $(z_0,x_0, y_{0})\in S^{3*}$ and $l$ is proper.

It remains to prove that $l$ is onto. First, 
for every $(x,y)$, 
$$
L_{(x,y)}=l(S^{3*})\cap M_{\dot\eta(x,y)},
$$
is a closed interval, being the image of an interval. If $(\gamma_+,\gamma_-)$ is a fixed point of $\gamma$, then $L_{(\gamma_+,\gamma_-)}$ is invariant by $\rho(\gamma)$ that acts as a translation on $M_{(\gamma_+,\gamma_-)}$. It follows that 
$$
L(\gamma_+,\gamma_-)=M_{(\gamma_+,\gamma_-)}.
$$ 
Since the set of fixed points of element of $\grf$ is dense in $\bgrf^2$, and $l(S^{3*})$ is closed by properness, we conclude that $l(S^{3*})=M$ and that $l$ is onto.
\qed

As a corollary, we may now prove Proposition \ref{mgcompact}
\vskip 0,3truecm
{\sc Proof of Proposition \ref{mgcompact}:}

Since $\grf$ acts properly on $S^{3*}$, and $l$ is proper and onto, the action of $\grf$ on $M$ in proper: indeed for any compact $K$ in $M$
$$
\{\gamma\in\grf \ \gamma(K)\cap K\not=\emptyset\}\subset\{\gamma\in\grf, \  \gamma(f^{-1}(K))\cap f^{-1}(K)\not=\emptyset\},
$$
Hence,
$$
\sharp\{\gamma\in\grf, \  \gamma(K)\cap K\not=\emptyset\}\leq\sharp\{\gamma\in\grf, \  \gamma(f^{-1}(K))\cap f^{-1}(K)\not=\emptyset\}<\infty.
$$
Next since $\grf$ has no torsion elements and acts properly, it follows that the action of $\grf$ on $M$ is free.
The space $M/\grf$ is therefore a topological manifold. Finally, the quotient map of $l$ from $S^{3*}/\grf$ (which is compact) being onto, it follows $M/\grf$ is compact.

\subsubsection{Contracting the leaves}

We notice that $M$ is topologically a trivial bundle.  Let $\mathcal L$ be the foliation of $M$ by parallel vector fields along the curves $c_s: t\mapsto(s,t)$ in $D_\kappa$. We choose a $\grf$-invariant metric on $M$. We prove now:

\begin{proposition} \label{contractleaves}
The action $\phi_t$ of $\mathbb R$ on $M$ contracts the leaves of $\mathcal L$.
\end{proposition}

\proof In the proof of Proposition \ref{mgcompact}, we exhibited a proper onto continuous $\grf$-equivariant map $l$  from $\bgrf^{3*}$ to $M$. This map is such that 
\begin{equation}
l(z,x,y)\in M_{\dot\eta(x,y)}.\label{defl}
\end{equation}

From  Proposition \ref{factnegcurv2}, by choosing a hyperbolic metric on $S$ we have a  flow $\psi_t$ by proper homeomorphisms on $\bgrf^{3*}$ such that
\begin{enumerate}
\item $\psi_{t}(z,x,y)=(w,x,y)$,
\item let $\seq{t}$ be a sequence of real numbers going to infinity, let $(z,x,y)\in\bgrf^{3*}$, let $\seq{\gamma}$ be a sequence of elements of $\grf$, such that 
$$
\{\gamma_{m}\circ\psi_{t_{m}}(z,x,y)\}\inn \hbox{ converges to } (z_0,x_{0},y_{0}).
$$ 
Then for any $w$, $v$ such that $(w,x,v)\in\bgrf^{3*}$, there exists $v_{0}$ such that
$$
\{\gamma_{m}\circ\psi_{t_{m}}(w,x,v)\}\inn\hbox{ converges to } (w_0,x_{0},y_{0}).
$$ 
\end{enumerate}

From the compactness of $S$, the properness of $l$ and the flows, and we obtain positive constants $a$ and $b$ such that
$$
\forall u, \forall T, \exists t\in]T/a-b,Ta+b[, l( \psi_t(u))=\phi_T(l(u)).
$$

Therefore, we deduce the following assertion holds: let $\seq{t}$ be a sequence of real numbers going to infinity, let $u\in M$, let $\seq{\gamma}$ be a sequence of elements of $\grf$, such that 
$$
\{\gamma_{m}\circ\phi_{t_{m}}(u)\}\inn \hbox{ converges to } u_0\in M_{\dot\eta(x_0,y_0)}.
$$ 
Then for any $w$, $s$ such that $w\in M_{\dot\eta (x,s)}$, there exists $v_{0}$ such that
$$
\{\gamma_{m}\circ\phi_{t_{m}}(w)\}\inn\hbox{ converges to } v_{0} \in M_{\dot\eta(x_0,y_0)}.
$$ 
Hence the hypothesis of Proposition \ref{conditionalternate} holds for the $\mathbb R$-action on $M$ and in particular this action contracts the leaves of $\mathcal  L$. \qed

\section{Appendix A : Filtrated Spaces, Holonomy.}\label{applam}

We are going to describe a notion of a topological space with "nested" laminations. This require some definitions.
First we introduce some notations,  Let 
$$
Z=Z_1\times\ldots\times Z_p
$$
For $k<p$, we note $p_k$ the projection
$$
Z\rightarrow  Z_1\times\ldots\times Z_k
$$
If $U$ is a subset of $Z$ and $x$ a point of $U$, we define the {\em $k^{\hbox{th}}$-leave}  through $x$ in $U$ to be
$$
U^k_x= U\cap p_k^{-1}\{p_k(x)\}.
$$
We notice that $U^{k+1}_x\subset U^k_x$. The {\em higher dimensional leaf} is $U_{1}$.
If $\phi$ is a map from $U$ to to a set $V$, we note
$$
\phi^k_x=\phi\vert_{U^k_x}.
$$
If $V\subset W=W_1\times\ldots\times W_k$, we say $\phi$ is a {\em filtrated map} if 
$$
\forall k,\ \phi(U^k_x)\subset V^k_{\phi(x)}.
$$
\subsection{Filtrated Space, Lamination}\label{aff}

In this section, we define a notion of filtrated space for which it does make sense to say some maps are "smooth along leaves with derivatives varying continuously".

\begin{definition}{\sc [$C^\infty$-filtrated space]} We say a metric space $P$ is {\em $C^\infty$-filtrated} if there exist 
\begin{itemize}
\item a covering of $P$ by open sets $U_i$, called {\em charts}, 
\item Hölder homeomorphisms, called {\em coordinates}, 
$\phi_i$ of $U_i$ with $V_1^i\times V_2^i\times\ldots\times V^i_p$ where for $k>1$, $V^k_i$ is an open set in a finite dimensional affine space, \end{itemize}
We furthermore make the following assumption about coordinates changes. Let 
$$
\phi_{ij}=\phi_i\circ(\phi_j)^{-1},
$$
defined from 
$
W_{ij}= \phi_j(U_i\cap U_j),
$ 
to $W_{ji}$.
We suppose
\begin{itemize}
\item $\phi_{ij}$ is a filtrated map,
\item $(\phi^1_{ij})_x$ is a $C^\infty$-map whose derivatives depends continuously on $x$. 
\end{itemize}
\end{definition}
We observe by the last assumption, for all $k$,  $(\phi^{k}_{ij})_x$ is a $C^\infty$-map whose derivatives depends continuously on $x$.

When $p=2$, we speak of a {\em laminated space}. 
We may want to specify the number of nested leaves, in which case we talk of $p$-filtrated objects. 

The Hölder hypothesis is somewhat irrelevant, but is needed in the applications we have in mind.

We now extend the definitions of the previous paragraph.
Let $P$ be a $C^\infty$-$p$-filtrated space. Let $k<p$.
\begin{definition}{\sc [Leaves] }
 The {\em $k^{\hbox{th}}$-leaves} are the equivalence classes of the equivalence relation generated by
$$
y \mathcal R x \iff p_k(\phi^i(y))=p_k(\phi^i(x)).
$$ 
\end{definition}
Let $P$ and $\overline P$ be two $C^\infty$-$p$-filtrated spaces.  Let $\phi_i$ be coordinates on the $P$ and $\overline\phi_j$ be coordinates on $\overline P$.

\begin{definition}{\sc [Filtrated maps and immersions]}
Let $\psi$ be a map from $\overline P$ to $P$. We say $\psi$ is  a {\em $C^\infty$-Filtrated maps} if
\begin{itemize}
\item $\psi$ is Hölder,
\item $\psi$ send $k^{\hbox{th}}$-leaves into $k^{\hbox{th}}$-leaves, 
\item $(\phi_i\circ\psi\circ\overline\phi_j)^k_x$ are $C^\infty$ and their derivatives vary continuously with $x$. 
\end{itemize}

A {\em filtrated immersion} is a filtrated map $\psi$ whose leafwise tangent map is injective: in other words, $(\phi_i\circ\psi\circ\overline\phi_j)^k_x$ is an immersion.
\end{definition}

\begin{definition}{\sc [Convergence of filtrated maps]}
Finally let $\{\psi_n\}_{n\in \mathbb N}$ be a sequence of $C^\infty$-filtrated maps from $\overline P$ to $P$. We say it {\em $C^\infty$ converges  on every compact set} if 
\begin{itemize}
\item it converges uniformly on every compact set
\item all the derivatives of $(\phi_i\circ\psi_n\circ\overline\phi_j)^k_x$  converges uniformly on every compact set (as a function of $x$).
\end{itemize}
\end{definition}

\begin{definition}{\sc [Affine Leaves]}
We say that $P$ is $C^\infty$-laminated by {\em affine leaves} or carries a {\em leafwise affine structure} if furthermore $\phi_{ij}\vert_{U_x}$ is an affine map.
\end{definition}

\subsection{Holonomy Theorem}
We now prove the following result which generalises Ehresmann-Thurston Holonomy Theorem.

\begin{theorem}\label{laminholo}
Let $P$ be a $p$-filtrated space. Let $G$ be a group of $C^\infty$-filtrated leafwise immersions of $P$, equipped with the topology of $C^\infty$-convergence on compact sets. Let $V$ be a compact $p$-filtrated space and  $\tilde V$ be a connected Galois covering with finitely generated Galois group $\Gamma$.  Let $U$ the set of homomorphisms $\rho$ from $\Gamma$ to $G$ such that there exists a $\rho$-equivariant filtrated immersion from $\tilde V$ to $P$. Then $U$ is open

Moreover, suppose $\rho_0$ belongs to $U$. Let $f_0$ be a $\rho_0$-equivariant filtrated immersion from $V$ to $P$. If $\rho$ is close enough to $\rho_0$, we may choose a  $\rho$-equivariant filtrated immersion $f$ arbitrarily close to $f_0$ on compact sets. 
\end{theorem}

We note again that we could have ripped the word ``Hölder"  from all the definitions and still have a valid theorem.

\proof  We choose a finite covering $U^i$ of $V$ such that
\begin{itemize}
\item $U^i$ are charts on $P$,
\item $U^i$ are trivialising open sets for the covering $\pi:\tilde V\rightarrow V$.
\end{itemize}
We choose  an open chart $U^1$ as well as a subset  $\tilde U^1\subset \tilde V$, such that $\pi$ is a homeomorphism from $\tilde U^1$ to $U^1$. We make the following temporary definitions.

A {\em loop} is a sequence of indices $i_1,\ldots i_l$ such that $i_l=i_1=1$ and
$$
U^{i_{j}}\cap U^{i_{j+1}}\not=\emptyset.
$$

A loop defines uniquely a sequence of open sets $\tilde U^{j}$ such that
\begin{itemize}
\item $\pi$ is a homeomorphism from $\tilde U^{j}$ to $U^{i_j}$,
\item $\tilde U^{{j}}\cap\tilde U^{{j+1}}\not=\emptyset $.
\end{itemize}

A loop $i_1,\ldots i_l$ is {\em trivialising} if $\tilde U^{1}=\tilde U^{l}$. In general, we associate to a loop  the element $\gamma$ of $\Gamma$ such that $\gamma(\tilde U^{1})=\tilde U^{p}$. The group $\Gamma$ is  the group of loops (with the product structure given by concatenation) modulo trivialising loops. This is just a way to choose a presentation of $\Gamma$ adapted to the charts $U^i$.

A {\em cocycle} is a finite sequence of $g=\{g^{ij}\}$ of elements of $G$ such that for every trivialising loop $i_1,\ldots, i_p$ we have
$$
g^{i_1i_2}\ldots g^{i_{p-1}i_p}=1.
$$
It follows that every cocycle $g$ defines uniquely a homomorphism $\rho_g$ from  $\Gamma$ to $G$, and furthermore the map $g\mapsto \rho_g$ is open. 

Let $g$ be a  cocycle, a $g$-equivariant map  $f$ is a finite collection $\{f_i\}$ such that 
\begin{itemize}
\item $f_i$ is a filtrated map from $U^i$ to $P$,
\item on $U^i\cap U^j$, $f_j=g^{ij}f_i$.
\end{itemize}
It is easy to check there is a one to one correspondence between $\rho_g$-equivariant laminated maps and $g$-equivariant maps.

We now oberve the following fact which follows from the existence of partition of unity:
\vskip 0.5truecm
\noindent{\em  Let $W_0,W_1, W_2$ be three open set in $V$, such that
$$
\overline{W_0}\subset W_1\subset \overline{W_1} \subset W_2.
$$ Let  $h$ a filtrated map defined on $W_2$, then there exists $\epsilon$ such that if $h_1$ is filtrated map defined on $W_1$ and $\epsilon$-close to $h$ on $W_1$ then there exists $h_0$, $2\epsilon$-close to $h$ on $W_2$, which coincides with  $h_1$ on $W_0$.}
\vskip 0.5truecm
Let us now begin the proof. Let $g$ be a cocycle  associated to a covering $\mathcal U=(U_1,\ldots,U_m)$. Let $f$ be a $g$-equivariant immersion. Let now $\overline{g}$ be a cocycle arbitrarily close to $g$.

We proceed by induction to build a $\overline{g}$-equivariant map  $\overline{f}$ defined on a smaller covering $\mathcal V= (V_1,\ldots,V_m)$ and close to $f$. Suppose
$$
\overline{f}^i=\{\overline{f}^i_1,\ldots,\overline{f}^i_{i-1}\}
$$ 
is a $\overline{g}$-equivariant map defined on 
$\mathcal V^i=(V^i_1,\ldots,V^i_{i-1}),$ with $V^i_k\subset U_k$.

Suppose also that $\overline f$ is close to $f$ on $V^i_l$, with $l<i$. Let
$$
W_i=V^i_1\cup\ldots\cup V^i_{i-1}.
$$
We now build $\overline{f}^i_i$ on a smaller subset $V^i_i$ of $U_i$ in the following way:
let $h_{i}=\overline{g}_{il}  \overline{f}^i_l$ on $W_i\cap U_{i}$.
The map $h_{i}$ is well defined on $W_i\cap U_i$ and close to $f_i$. We use our preliminary obervation to build  $\overline{f}^i_{i}$ close to $f_i$ on  $U_i$ and coinciding with
$h_i$ on a slightly smaller open subset $Z_i$ of $W_i\cap U_i$. We finally define
$V^{i+1}_l=V^i_l\cap Z_i$. This completes the induction.

In the end, we obtain a $\overline{g}$-equivariant map $\overline{f}$ defined on slightly smaller open subsets of $U_i$, and close to $f$. Therefore, it follows $\overline{f}$ is an immersion.

The construction above proves also the last part of the statement.
\qed

\subsection{Completeness of affine structure along leaves}

Let  $V$ be a compact space $C^\infty$-laminated by affine leaves. It follows that every leaf carries an affine structure. We shall say $V$ is {\em leafwise complete} if
the universal cover of every leaf is isomorphic, in the affine category, to the affine space.

We want to prove the following

\begin{lemma}\label{affcomp}
Let $V$ be a compact space  $C^\infty$-laminated by affine leaves. Let $E$ be the vector bundle over $V$ whose fibre at $x$ is the tangent space at $x$ of  the leaf $\mathcal L_x$.
Assume there exists a one-parameter group $\phi_t$ of homeomorphisms of the leaves such that
\begin{itemize}
\item for every leaf $\mathcal L_x$, $\phi_t$ preserves $\mathcal L_x$ and acts as a one parameter group of translation on  $\mathcal L_x$, generated by the vector field $X$.
\item Let $L=\mathbb R.X$, then the action of the lift of $\phi_t$ on the vector bundle 
$$
F=E/L
$$ 
is uniformly contracting.
\end{itemize}
Then $V$ is leafwise complete.
\end{lemma}

\proof  For every $x$ in $V$, 
$$
O_x=\{u\in E_x/ x+u\in \mathcal L_x.
$$
Let $$
O=\cup_{x\in V}O_x.
$$
We oberve that $O$ is an  open subset of $E$, which is invariant by $\phi_t$. By hypothesis, we have
$$
L\subset O.
$$
Since $\phi_t$ is contracting on $F=E/L$ and $V$ is compact, $L$ admits a $\phi=\phi_1$ invariant supplementary $F_0$. Let us recall the classical and well known proof of this fact. We choose a supplementary $F_1$ to $L$. Then $\phi^* F_1$ is the graph of an element $\omega$ in  $K=F_1^*\otimes L$.  We now identify $F_1$ with $F$ using the projection. Since the action of $\phi_t$ is uniformly contracting on $F$, the action of $\phi_{-t}$ is uniformly contracting on $K=F_1^*\otimes L$, hence exponentially contracting by compactness. It follows the following element of $K=F_1^*\otimes L$ is well defined:
$$
\alpha= \sum_{p=-1}^{p=-\infty}(\phi^p)^*\omega.
$$
This section $\alpha$ satisfies the {\em cohomological equation}
$$
\phi^*\alpha -\alpha= \omega. 
$$
This last equation exactly means that the graph $F_0$ of $\alpha$ is $\phi$ invariant.

Let $u\in E$. Write $u=v+\lambda X$, with $v\in F_0$. It follows that for $n$ large enough, we have $\phi^n (u)\in O$. Since $O$ is invariant by  $\phi$, we deduce that $u\in O$. Hence $O=E$ and $V$ is leafwise complete. \qed 

\section{Appendix B : the symplectic nature of cross ratio}\label{sympnat}
In this appendix, we explain how to construct different cross ratio from hyperconvex curves from a symplectic point of view. We concentrate on the case of projective spaces, although the construction can be extended to flag manifolds to produce a whole family of cross ratios. However in this case, for the moment, we do not know how to characterise these cross ratios using functional relations, as we did in the case of curves in the projective space. 

We also give more precisions concerning cross ratios associated  of hyperconvex curves and prove the result used in the proof.

\subsection{A symplectic construction}
All the examples of cross ratio we have defined may be interpreted from the following  "symplectic" construction. Let $V$ and $W$ be two manifolds of the same dimension and $O$ be an open set of $V\times W$ equipped with an {\em exact} symplectic structure, or more generally an exact  two-form $\omega$. We assume furthermore that the two foliations coming from the product structure
\begin{eqnarray*}
\mathcal F^+_w&=&O\cap(V\times\{w\})\\
\mathcal F^-_v&=&O\cap (\{v\}\times W),
\end{eqnarray*}
satisfy the following properties
\begin{enumerate}
\item Leaves are connected.
\item The first cohomology groups of the leaves are reduced to zero.
\item $\omega$ retricted to the leaves is zero.
\item finally, let {\em squares} be closed curves of the form $c_1\cup c_2\cup c_3\cup c_4$ where $c_1$ and $c_3$ are along $\mathcal F^+$ and $c_2$ and $c_4$ are along $\mathcal F^-$; we assume that squares are homotopic to zero.
\end{enumerate}

\vskip 1 truecm
\noindent\rmks
\begin{itemize}\label{afflagr}
\item When $\omega$ is symplectic, by definition every leaf $\mathcal F^\pm$ is lagrangian. Then it is a standard fact that it carries a flat affine structure. If every leaf is simply connected and complete from the affine point of view, condition (4) above is satisfied: one may "straighten" the edges of the square, that is deforming them into geodesics of the affine structure, then use these straightening to define a homotopy.
\item We shall give later on  examples of this situation.
\end{itemize}

Associated to the above data is a function $B$, called the {\em polarised cross ratio} defined on
$$
U=\{(e,u,f,v)\in V\times W\times V\times W/ (e,u),(f,u),(e,v),(f,v)\in O\},
$$
in the following way. We consider a map $G$ from the square $[0,1]^2$ to $O$ such 
that
\begin{itemize}
\item the image of the vertexes $(0,0)$, $(0,1)$, $(1,1)$, $(1,0)$ are respectively
$(e,u)$, $(f,u)$, $(e,v)$, $(f,v)$,
\item the image of every edge on the boundary of the square lies in a leaf of $\mathcal F^+$ or $\mathcal F^-$.
\end{itemize}
We define  the {\em polarised cross ratio} to be the function defined on  $U$ by  
$$
B(e,u,f,v)=e^{\frac{1}{2}\int_G\omega}.
$$
It is easy to check that the definition of $B$ does not depend on the choice of the specific map $G$.

Let now $\xi$ and $\xi^*$ be two maps from $S$ to $V$ and $W$ respectively such that for all distinct $x$ and $y$, $(\xi(x),\xi^*(y))$ lies in $O$. Then, we have the following immediate 
\begin{proposition}\label{sympcr}
The function $b$ defined by
$$
b(x,y,z,t)=B(\xi(x)),\xi^*(y),\xi(z),\xi^*(t)).
$$
satisfies
\begin{eqnarray*}
b(x,y,z,t)&=&b(z,t,x,y) \\
b(x,y,z,t)&=&{b(x,y,z,w)}{b(x,w,z,t)}\\
b(x,y,z,t)&=&b(x,y,w,t)b(w,y,z,t)\\
\end{eqnarray*}
\end{proposition}
This function $b$ is not defined for $x=y$ and $z=t$. It follows from the above proposition that  it extends to a crossratio provided that
$$
\lim_{y\rightarrow x}b(x,y,z,t)=0.
$$

We explain quickly a similar construction for triple ratio. We consider a sextuple 
$(e,u,f,v,g,w)$ in $V\times W\times V\times W\times V\times W$. Let now $\phi$ be a map from the interior of the regular hexagon $H$  in $V\times W$ such that the image of the edges lies in $\mathcal F^+$ or $\mathcal F^-$, and the  (ordered) image of the vertices are 
$(e,u)$, $(f,u)$, $(f,v)$, $(g,v)$, $(g,w)$, $(e,w)$. We check that the following quantity does not depend on the choice of $\phi$ :
$$
T(e,u,f,v,g,w)=e^{\frac{1}{2}\int_H\phi^*\omega}.
$$
Finally using the same notations as above, we check that
$$
t(x,y,z)=T(\xi(x),\xi^*(z),\xi(y),\xi^*(x),\xi(z),\xi^*(y)),
$$
is the triple ratio as defined in Paragraph \ref{triratio}.
\subsubsection{Period and action difference}

Let $\gamma$ be an exact Hamiltonian diffeomorphism of $O$. Let $\alpha$ and $\beta$ be two fixed points of $\gamma$ and $c$ a curve joining $\alpha$ and $\beta$. Since $\gamma$ is isotopic to the identity, it follows that $c\cup \gamma(c)$ bounds a disc $D$.  We define the {\em action difference} to be 
 $$
 \delta_\gamma(\alpha,\beta,c)=exp(\int_D \omega).
 $$ 
 We first recall the 
 \begin{proposition}
 The quantity $\delta:=\delta_\gamma(\alpha,\beta,c)$ just depends on the homotopy class of $c$.
\end{proposition}
\proof In our case, this follows from the fact  $\omega$ is exact. In general the action difference depends also on a path joining $\gamma$ to the identity. \qed

In our case, we have a preferred class of curves joining two points as we now explain: let  $\alpha=(a,b)$ and $\beta=(\bar a,\bar b)$ be two points of $O$. We notice that since squares are homotopic to zero, we have a well defined homotopy class $c_{a,b,\bar a,\bar b}$ for curves from $(a,b)$ to $(\bar a,\bar b)$, namely curves homotopic to $c^+\cup c^-\cup \bar c^{+}$, where $c^+$ is a curve along $\mathcal F^+$ going from $(a,b)$ to  $(y,b)$,  $c^-$ a curve along   $\mathcal F^-$ going from $(y,b)$ to  $(y,\bar b)$, and $c^-$ a curve along   $\mathcal F^+$ going from $(y,\bar b)$ to  $(\bar a,\bar b)$. 
By convention we set
$$
\delta_\gamma(\alpha,\beta)=\delta_\gamma(\alpha,\beta,c_{a,b,\bar a,\bar b}).
$$
Let $\phi=(\rho,\bar\rho)$ be a representation of $\grf$ in the group of Hamiltonian diffeomorphism of $O$ which are restriction of elements of $Di\!f\!\!f(V)\times Di\!f\!\!f(W)$. Let $(\xi,\xi^*)$ be a $\phi$-equivariant map of $\bgrf$ in $O$.   We now prove using the notations of the previous paragraph
\begin{proposition}
Let $\gamma$ be an element of $\pi_{1}(S)$ then,
$$
b_{\xi,\xi^*}(\gamma^+,y,\gamma^-,\gamma y)^2=\delta_{\phi(\gamma)}\big((\xi(\gamma^+),\xi^*(\gamma^-)),(\xi(\gamma^-),\xi^*(\gamma^+))\big).
$$
In particular,
$$
l_{b_{\xi,\xi^*}}(\gamma)=\frac{1}{2}\log\vert \delta_{\rho(\gamma)}\big((\xi(\gamma^+),\xi^*(\gamma^-)),(\xi(\gamma^-),\xi^*(\gamma^+))\big)\vert.
$$
\end{proposition}
\proof Let $f=(g,\bar g)$ be a Hamiltonian diffeomorphism of $O$, restriction of an element of $Di\!f\!\!f(V)\times Di\!f\!\!f(W)$. Let $(a,b)$ and $(\bar a, \bar b)$ be two fixed points of $f$. Let  as before $c=c^+\cup c^-\cup \bar c^{+}$ composition of
\begin{itemize}
\item
 $c^+$  a curve along $\mathcal F^+$  from $(a,b)$ to  $(y,b)$, 
 \item $c^-$ a curve along   $\mathcal F^-$ from $(y,b)$ to  $(y,\bar b)$,
 \item and $c^-$ a curve along   $\mathcal F^+$ from $(y,\bar b)$ to  $(\bar a,\bar b)$.
 \end{itemize}Assume $(a,b)$ and $(\bar a \bar b)$ be fixed points of $f$. Then $c\cup\gamma(c)$ is a ``square", {\it i.e} the composition of
\begin{itemize}
\item  a curve along $\mathcal F^+$ from $f(y,b)=(g(y), b)$ to $(y,b)$,
\item  a curve along $\mathcal F^-$ from $(y,b)$ to $(y,\bar b)$,
\item  a curve along $\mathcal F^+$ from $(y,\bar b)$ to $f(y,\bar b)=(g(y),\bar b)$ ,
\item  a curve along $\mathcal F^+$ from $(g(y),\bar b)$ to $(g(y), b)$
\end{itemize}

Let $D$ be a disk whose boundary is this square. By definition
$$
\delta_f\big((a,b),(\bar a,\bar b)\big)=\int_D\omega.
$$
The proposition follows from the definition of the cross ratio associated to $(\xi,\xi^*)$ when we take $f=(\rho(\gamma),\rho^*(\gamma))$ and
$$
(a,b,\bar a,\bar b)=(\xi(\gamma^+),\xi^*(\gamma^-),\xi(\gamma^+),\xi^*(\gamma^-)).
$$\qed 

\subsection{Projective spaces}
As a specific example of the previous situation, we wish to discuss the following case which makes the link with Section \ref{crosscurve}. Let $E$ be a vector space. We identify $\pnd$ with the set of hyperplanes of $E$. Let
$$
{\mathbb P}^{2*}=\pn\times\pnd\setminus\{(D,P)/D\subset P\}
$$
Using the identification of  $T_{(D,P)}{\mathbb P}^{2*}$ with 
$\Hom(D,P)\oplus \Hom(P,D)$, let
$$
\Omega((f,g),(h,j))=tr(f\circ j)-tr(h\circ g).
$$
Let $L$ be the $\mathbb R$-bundle over ${\mathbb P}(n)^{2*}$, whose fibre at $(D,P)$ is
$$
L_{(D,P)}=\{u\in D, f \in P^{\perp}/ \langle f, u\rangle=1\}/\{+1,-1\}.
$$
Then 
\begin{proposition}\label{sympproj}
There exists a connection form $\beta$ on $L$ such that
\begin{itemize}
\item Its curvature is symplectic and equal to $\Omega$.
\item Let $u \in D\subset E$, then the section 
$$
\xi_u : P \mapsto (u,f) \hbox{ such that } \langle u,f\rangle =1, \hbox{ and }  f\in P^\perp
$$
is parallel for $\beta$ above $\{D\}\times(\pnd \setminus\{P/ D\subset P\})$. 
\item If $h\in \mathbb P (F)$, we denote by $\hat h$ a nonzero  element of $h$. The polarised  cross ratio associated to $2\Omega$  is 
$$
b(u,f,v,g)=\frac{\langle\hat f ,\hat v\rangle\langle\hat g,\hat u\rangle}{\langle\hat f,\hat u\rangle\langle\hat g,\hat u\rangle}.
$$
\item  Let $f$ be an element of $SL(n,\mathbb R)$. Let $D$ (resp $\bar D$) be an eigenspace of dimension one for the eigenvalue $\lambda$ (resp. $\mu$). Then $(D,\bar D^\perp)$ is a fixed point of $f$ in $\mathbb P^{2*}$. The action of $f$ on $L_{(D,\bar D^\perp)}$ is the translation by $$\log\vert\lambda/\mu\vert.$$
\end{itemize}
\end{proposition}
\proof We consider the standard symplectic form $\Omega^0$ on $E\times E^*/\{+1,-1\}$. We observe that $\Omega^0=d\beta^0$ where $\beta^0_{(u,f)} (v,g)=\langle u,g\rangle$.
We have a symplectic action of $\mathbb R$ given by
$$
\lambda.(u,f)=(\lambda^{-1} u, \lambda f),
$$
with moment map
$$
\mu\big((u,f)\big)=\langle  f, u\rangle.
$$ 
We observe that $L=\mu^{-1} (1)$. Therefore, we obtain that $\beta=\beta_0\vert_L$ is a connection form for the $\mathbb R$-action, whose curvature $\Omega$ is the symplectic form obtained by reduction of the Hamiltonian action of $\mathbb R$.

We now compute explicitely $\Omega$. Let $(D,P)$ be an element of $\mathbb P^{2*}$. Let $\pi$ be the projection onto $P$ parrallel to $D$.
We  identify  $T_{(D,P)}{\mathbb P}^{2*}$ with $\Hom(D,P)\oplus \Hom(P^\perp,D^\perp)$. Let $(f,\hat g)$ be an element of $T_{(D,P)}{\mathbb P}^{2*}$. Let $u\in D$, $\alpha\in P^\perp$ with $(u,\alpha)\in L$, then $(f(u),\hat g(\alpha))$ is an element of $T_{(u,\alpha)}L$ which projects to $(f,\hat g)$.
By definition of the symplectic reduction if $(f,\hat g)$ and $(h,\hat l)$ are elements of  $T_{(D,P)}{\mathbb P}^{2*}$, then
$$
\Omega((f,\hat g),(h,\hat l))=\langle \hat l(\alpha),f(u)\rangle -\langle \hat g(\alpha),h(u)\rangle.$$
Finally, let $\pi$ be the projection onto $P$ in the $D$ direction. We define
$$
\mapping{\Hom(P,D)}{\Hom(P^\perp,D^\perp)}{f}{\hat f=(f\circ \pi )^*}.
$$
In particular
$$
\langle \hat l(\alpha),f(u)\rangle =tr(l\circ f) \langle \alpha,u\rangle
$$
The second point follows immediately by the explicit formula for $\beta^0$. Using this, a computation of the holonomy of this connection shows the formula about the cross ratio.

The last point is obvious.
\qed

\auteur
\end{document}